\documentclass[11pt]{article}

\usepackage{natbib}
\usepackage[english]{babel}
\usepackage[utf8]{inputenc}
\usepackage{theorem}
\usepackage{amsmath,amssymb}
\usepackage{amsfonts,graphicx,verbatim,psfrag, graphics}
\usepackage{verbatim, rotating, epsfig}
\usepackage{mathtools}
\usepackage{bbm, multicol}
\usepackage[dvipsnames]{xcolor}
\usepackage[hidelinks]{hyperref} % make internal reference to bib
\usepackage[T1]{fontenc}
\usepackage{lmodern}
\usepackage[normalem]{ulem}

\DeclareMathAlphabet{\mathpzc}{OT1}{pzc}{m}{it}

\newcommand{\be}{\begin{equation}} \newcommand{\ee}{\end{equation}}
\newcommand{\bea}{\begin{eqnarray}} \newcommand{\eea}{\end{eqnarray}}

\newcommand{\ba}{\begin{array}} \newcommand{\ea}{\end{array}}

\newcommand{\noi}{\noindent}

\newcommand{\Bthe}{\begin{theorem}}
\newcommand{\Ethe}{\end{theorem}}

{\theoremstyle{plain}} \theoremheaderfont{\scshape}

\newtheorem{theorem}{Theorem}

\newtheorem{corollary}{Corollary}
\newtheorem{remark}{Remark}
\newtheorem{proposition}{Proposition}
\newtheorem{lemma}{Lemma}

\newtheorem{assumption}{Assumption}

\newtheorem{definition}{Definition}

\newenvironment{proof}{\noi {\it Proof.} }

\newtheorem{assumptionalt}{Assumption}[assumption]

\newenvironment{assumptionp}[1]{
  \renewcommand\theassumptionalt{#1}
  \assumptionalt
}{\endassumptionalt}

\begin{document}

\title{\centering { \scshape Existence, uniqueness, and regularity of\\ solutions to nonlinear and non-smooth parabolic obstacle problems\thanks{Emails: theod@illinois.edu and b-strulovici@northwestern.edu. \\
We thank Jose Blanchet, Benjamin Bernard, Svetlana Boyarchenko, Simone Cerreia Vioglio, Ibrahim Ekren, R. Vijay Krishna, Hugo Lavenant, Massimo Marinacci, Dylan Possaima\"{i}, Mete Soner, Mehdi Talbi, Nizar Touzi, St\'{e}phane Villeneuve, and seminar participants at the the Mathematical Finance seminar at Florida State University, the ``Applications of Stochastic Control to Finance and Economics'' workshop hosted by the Banff International Research Station for Mathematical Innovation and Discovery, and at Bocconi University (Decision Science).}}}
\author{Th\'{e}o Durandard and Bruno Strulovici \\ UIUC and Northwestern University}
\date{\today}

\maketitle

\begin{abstract}
    We establish the existence, uniqueness, and $W^{1,2,p}$-regularity of solutions to fully-nonlinear, parabolic obstacle problems when the obstacle is the pointwise supremum of functions in $W^{1,2,p}$ and the nonlinear operator is required only to be measurable in the state and time variables. In particular, the results hold for all convex obstacles. Applied to stopping problems, they provide general conditions under which a  decision maker never stops at a convex kink of the stopping payoff. The proof relies on new $W^{1,2,p}$-estimates for obstacle problems when the obstacle is the maximum of finitely many functions in $W^{1,2,p}$.
\end{abstract}

\maketitle
%%%%%%%%%%%%%%%%%%%%%%%%%%%%%%%%%%%%%%%%%%%%%%%%%%%%%%%%%%%%%%%%%%%%%%

% Text of your paper here

\section{Introduction}\label{sec:Introduction}

We study a fully nonlinear parabolic obstacle problem with Dirichlet boundary data on the domain $\mathcal{Y} = [0,T) \times \mathcal{X}$  where $T$ is finite and $\mathcal{X}$ is a bounded, open subset of $\mathbb{R}^d$ for some $d \in \mathbb{N}$: 
\begin{equation}\label{eq:obstacle}
    \begin{cases}
        \max \left\{ u_t + F\left(t,x,u, u_x, u_{xx}\right), \,g - u  \right\} = 0 \text{ on } \mathcal{Y}, \\
        u=b \text{ on } \partial \mathcal{Y}.
    \end{cases}
\end{equation}
Here, $\partial \mathcal{Y}= \{T\} \times \mathcal{X} \; \cup [0,T] \times \partial \mathcal{X}$ is the boundary of $\mathcal{Y}$ and $F$ is a measurable nonlinear uniformly elliptic operator defined on $\mathcal{Y}\times \mathbb{R}  \times \mathbb{R}^{d} \times \mathbb{S}^d$.\footnote{$\mathbb{S}^d(\mathbb{R})$ denotes the set of symmetric $d \times d$ real-valued matrices.} 

{Parabolic obstacle problems} consist of finding the smallest function that (i) exceeds a given obstacle function and (ii) is a supersolution of a given parabolic equation \citep{petrosyan2007parabolic}. Such a function solves an equation of the form \eqref{eq:obstacle} and is called the solution of the obstacle problem. Obstacle problems arise in physics, e.g., to study phase transitions (\cite{lame1831memoire,stefan1889einige,friedman1982variational}), in biology, e.g., to study tumor growth (\cite{greenspan1976growth,bazaliy2003free}), in economics, e.g., to study learning and investment decisions (\cite{wald2004sequential,dixit1993art,moscarini2001optimal,decamps2024investment}), and in finance, e.g., to study American options (\cite{jacka1991optimal,shiryaev1999essentials,villeneuve1999exercise}). 

We are particularly interested in obstacle problems that arise in  optimal stochastic control and stopping problems. When the primitives of the corresponding obstacle problem are sufficiently regular, it is well-known that the value function for the optimal control and stopping of a diffusion is the solution of an obstacle problem where the operator $F$ is the Hamilton-Jacobi-Bellman operator and the obstacle is the stopping payoff \citep{bensoussan1978inegalites,friedman1982variational, karatzas2001controller,peskir2006optimal,petrosyan2012regularity,strulovici2015smoothness}. In many applications, however, the operator, the obstacle, or the domain that define the obstacle problem fail to be regular.

For example, the obstacle $g$ that arises in stopping problems often takes the form $g = \underset{a \in A}{\sup } \, g^a$, and may thus have kinks even when the individual functions $g^a$ are smooth. This form arises when the decision maker must decide, upon stopping, on some action $a$. In this case, $g^a(t,x)$ is the stopping payoff when stopping at time $t$ in state $x$ and taking action $a$ upon stopping. A decision maker who stops at time $t$ would then optimally choose the action $a$ that maximizes $g^a(t,x)$ over all possible actions $a\in A$, and the resulting stopping payoff, i.e., obstacle, is $g$. Such problems are pervasive in information acquisition models where  ``stopping'' means concluding the information acquisition stage and $a$ is the decision taken after that stage (\cite{wald1992sequential}, \cite{decamps2006irreversible}, \cite{fudenberg2018speed} and \cite{camboni2023under}). Similarly, in many obstacle problems that arise in economics, the domain of the state variable fails to be smooth. In information acquisition problems, for example, the state variable is a belief and the relevant domain is a probability simplex, which has kinks at its vertices (and, more generally, along the lower-dimensional faces of the probability simplex). As another example, the state variable in many finance and general equilibrium problems is a price vector whose relevant domain is a positive orthant, which again fails to be smooth. 

To address the primitives' lack of regularity, researchers have often used techniques tailored to their specific application. For instance, \cite{decamps2006irreversible} consider the elliptic obstacle equation associated with an optimal stopping problem when the payoff upon stopping (i.e., the obstacle) is the maximum of two smooth convex functions, and use a probabilistic local time argument to show that the decision-maker never stops at the convex kink. The results for regular obstacles then apply, and the authors can show that the value function is a smooth solution of \eqref{eq:obstacle}.

This paper aims to address these challenges all at once. Our main result, Theorem \ref{theorem:existenceobstacleproblemirregulardomain}, establishes the existence, uniqueness, and $W^{1,2,p}$-regularity of a solution to fully nonlinear obstacle problems when (i) the operator is required only to be {measurable} in $(t,x)$, (ii) the obstacle is the {supremum} of functions in $W^{1,2,p}$---and may thus have kinks---and (iii) the domain is required only to be  Lipschitz---and may thus have corners and other kinks. Theorem~\ref{theorem:existenceobstacleproblemirregulardomain} covers various applications considered in economics. For instance, it allows us to recover \cite{decamps2006irreversible}'s result when combined with a Sobolev embedding theorem.\footnote{Corollary \ref{corollary:smoothness} makes this point formally.} More generally, it guarantees that the solution of \eqref{eq:obstacle} is continuously differentiable in space and, thus, that the decision-maker never stops at a point of nondifferentiability of the stopping payoff $g$. Since the results apply to any convex obstacle, they imply that  any point of nondifferentiability of the stopping payoff must be in the continuation region.

Results similar to Theorem \ref{theorem:existenceobstacleproblemirregulardomain} have appeared under stronger assumptions. A precursor is \cite{friedman1982variational}, which shows that when the primitives are smooth and the operator $F$ is linear, the obstacle problem has a unique solution in $W^{1,2,\infty}\left( {\mathcal{Y}} \right)$. \cite{petrosyan2007parabolic,figalli2014general,indrei2016regularity,audrito2023regularity} extend this result to nonlinear operators when the obstacle and operator are sufficiently regular (the operator must be at least continuous and the obstacle must belong to $W^{1,2,\infty}$) and the operator satisfies strong convexity and growth conditions. When the operator is not smooth, solvability and regularity results for elliptic and parabolic equations in the absence of obstacles appear in \cite{caffarelli1989interior}, \cite{escauriaza1993w}, \cite{caffarelli1996viscosity}, \cite{winter2009w}, and \cite{krylov2010bellman} for the elliptic case, and in \cite{crandall1998remarks}, \cite{crandall1999existence}, \cite{crandall2000lp}, \cite{dong2013fully}, and \cite{krylov2017existence} for the parabolic case, among others. The recent monograph \cite{krylov2018sobolev} offers an up-to-date general treatment. Finally, \cite{byun2018nondivergence,byun2022w} study the regularity of viscosity solutions for obstacle problems when the obstacle is smooth and one of the following conditions hold: the operator is linear or the problem is elliptic.

Specifically, our main theorem generalizes existing results in two directions: it extends the best known regularity for nonlinear parabolic equations to nonlinear parabolic problems and it relaxes the regularity needed from the obstacle. \cite{crandall2000lp} and \cite{krylov2017existence,krylov2018sobolev} prove the existence, uniqueness, and $W^{1,2,p}$-regularity of solutions of fully nonlinear parabolic problems {in the absence of an obstacle} (let alone an irregular one). In addition, \cite{krylov2017existence,krylov2018sobolev} assume that the domain is smooth and \cite{crandall2000lp} assume that the operator is Lipschitz continuous in the derivatives of the value function, while the theorem of the present paper requires only continuity.\footnote{\cite{krylov2017existence,krylov2018sobolev} use the same concept of solutions as we do: $L^p$-solutions, while \cite{crandall2000lp} work with the weaker concept of $L^p$-{viscosity} solutions. However, these concepts coincide under our assumptions on the nonlinear operator.} \cite{byun2018nondivergence} proposes a new method to prove that parabolic obstacle problems on smooth domains are solvable in $W^{1,2,p}$ when the operator is linear and the obstacle is in $W^{1,2,p}$, a condition that is not satisfied by the kinked obstacles that arise in various stopping problems. Finally, \cite{byun2022w} adapt the technique in \cite{byun2018nondivergence} to study fully nonlinear {elliptic} obstacle problems in smooth domains when the obstacle is in $W^{2,p}$. 

Theorem \ref{theorem:existenceobstacleproblemirregulardomain} builds on the results and ideas in these earlier papers to establish the existence, uniqueness, and $W^{1,2,p}$-regularity of solutions to fully nonlinear parabolic obstacle problems when the operator is {measurable} in $(t,x)$, the obstacle is {irregular}, and the domain is not required to be smooth. The proof of Theorem \ref{theorem:existenceobstacleproblemirregulardomain} relies on PDE methods and generalizes existing results derived by probabilistic methods in the context of optimal stopping. Specifically, we extend the approximation argument developed in \cite{byun2018nondivergence} and \cite{byun2022w} to the fully nonlinear parabolic case and we allow for non-smooth obstacles. To do so, a key ingredient in our analysis is the Sobolev-space approach developed in \cite{krylov2018sobolev}, which provides $W^{1,2,p}$-regularity results for fully nonlinear parabolic equations under the minimal regularity assumptions on the operator we maintain. In particular, we first extend the $W^{1,2,p}$ a priori estimates for nonlinear parabolic Dirichlet problem {without obstacles} derived in \cite{krylov2018sobolev} (Theorem 12.1.7 and 15.1.3) to obstacle problems with regular $W^{1,2,p}$ obstacles and smooth domain by standard methods. We then show via an induction argument that the existence and regularity results for a single $W^{1,2,p}$ obstacle extend to the pointwise maximum of finitely many $W^{1,2,p}$ obstacles (Lemma \ref{theorem:existenceobstacleproblemfiniteA}). Finally, we extend the result to Lipschitz domains and to obstacles given by the supremum of a separable family of functions in $W^{1,2,p}$ by a limit argument. So, our contribution relative to \cite{krylov2018sobolev} is one of scope. First, we extend his regularity results for nonlinear parabolic equations to regular obstacle problems. Then, we establish the existence and regularity of solutions to fully nonlinear obstacle problems with irregular features -- most notably, kinked obstacles and non-smooth domains, which are essential in applications such as the stopping problem studied in Section \ref{sec:application} -- using a novel induction argument.

\subsection{Outline of the Paper}\label{subsec:outline}

Section \ref{sec:preliminaries} introduces notation and definitions. Section \ref{sec:mainresult} formally states our assumptions and main result, Theorem \ref{theorem:existenceobstacleproblemirregulardomain}. In Section \ref{sec:application}, we present an application to robust investment using the results developed in this paper. Section \ref{sec:proof} contains the proof of Theorem \ref{theorem:existenceobstacleproblemirregulardomain}. We begin by establishing two key auxiliary results in Section \ref{sec:comparisonandstability}: (i) a comparison principle for fully nonlinear parabolic obstacle problems with measurable ingredients (Proposition \ref{prop:comparison}), and (ii) a stability theorem for solutions in $W^{1,2,p}(\mathcal{Y})$ (Theorem \ref{theorem:continuityequation}). Next, Section \ref{sec:finiteA} proves our main result in a simpler setting in which the obstacle is the maximum of finitely many functions in $W^{1,2,p}$ and the domain $\mathcal{X}$ is smooth. Finally, Section \ref{sec:proofirregulardomains} establishes Theorem \ref{theorem:existenceobstacleproblemirregulardomain} in full generality by first extending the result of Section \ref{sec:finiteA} to general separable sets $A$, and then to non-smooth domains.

\section{Preliminaries}\label{sec:preliminaries}

\subsection{Notation}\label{subsec:notations}

\begin{itemize}
    \item For $D>0$, $B_D(x)$ denotes the open ball of radius $D$ centered around $x$ and $C_D(t,x)= [t,t+D) \times B_D(x)$.

    \item $\bar{\mathcal{O}}$ denotes the closure of $\mathcal{O}$ for the relevant topology.

    \item $d(\cdot, \cdot)$ is the Euclidean distance. For any sets $\mathcal{Y}, \mathcal{Y}' \subset \mathbb{R}^{d+1}$, define $diam\left(\mathcal{Y}\right)$ \,  $ = \sup \big\{ d\left((t,x), (t',x')\right) \, : \, (t,x), \, (t',x') \in \mathcal{Y}   \big\}$, and $dist\left(\mathcal{Y}, \mathcal{Y}'\right)= \inf\left\{ d\left((t,x), (t',x')\right) \, : \, (t,x) \in \mathcal{Y} \text{ and } (t',x') \in \mathcal{Y}'  \right\}$.

    \item $\mathcal{C}^k\left( \mathcal{Z}\right)$ is the space of continuous functions on $\mathcal{Z}$, if $k=0$, and of $k$-times continuously differentiable functions on $\mathcal{Z}$, if $k\geq 1$.

    \item $W^{1,2,p}\left(\mathcal{Z}\right)$ denotes the Sobolev space of functions defined on the set $\mathcal{Z}$ whose first weak time derivative and second weak space derivatives are $L^p$-integrable. $W^{1,2,p}_{loc}\left(\mathcal{Z}\right)$ is the space of functions that belong to $W^{1,2,p}\left(\mathcal{Z}'\right)$ for all compact subsets $\mathcal{Z}'$ of $\mathcal{Z}$.

    \item For a function $u\in W^{1,2,p}_{loc}\left(\mathcal{Z}\right)$, $u_t, u_x$, and $u_{xx}$ stand for the first weak derivatives of $u$ with respect to $t$, the weak gradient of $u$ with respect to $x$, and the weak Hessian of $u$ with respect to $x$, respectively.

    \item $\rightharpoonup$ denotes weak convergence.
    
    \item For $k\in \mathbb{N}$, an open bounded subset $\mathcal{Z}$ of $\mathbb{R}^d$ is $\mathcal{C}^{k, Lip}$ if there exists a map $\phi : \mathbb{R}^{d-1} \to \mathbb{R}$ that is Lipschitz continuous, if $k=0$, or $k$-times continuously differentiable with Lipschitz derivatives, if $k\geq 1$, such that, for every $x \in \partial \mathcal{Z}$, there exists a neighborhood $V$ of $x$, $\alpha >0, \beta >0$, an affine map $T: \mathbb{R}^{d-1} \times \mathbb{R} \to \mathbb{R}^d$ such that
	\begin{equation*}
	\partial {\mathcal{Z}} \cap V = T \left(\left\{ (\xi, \eta ) \in \mathbb{R}^{d-1} \times \mathbb{R} \, :  \, \xi \in B_{\alpha}(0) \text{ and } \eta = \phi(\xi) \right\} \right),
	\end{equation*}
	and
	\begin{equation*}
        \begin{aligned}
            & T \left(\left\{ (\xi, \eta ) \in \mathbb{R}^{d-1} \times \mathbb{R} \, :  \, \xi \in B_{\alpha}(0) \text{ and } \phi(\xi) < \eta < \phi(\xi) + \beta \right\} \right) \subseteq {\mathcal{Z}} \cap V \\
	\text{ and } & T \left(\left\{ (\xi, \eta ) \in \mathbb{R}^{d-1} \times \mathbb{R} \, :  \, \xi \in B_{\alpha}(0) \text{ and } \phi(\xi) -\beta < \eta < \phi(\xi) \right\} \right) \subseteq \mathbb{R}^d \setminus \left(\bar{{\mathcal{Z}}} \cap \bar{V}\right).
        \end{aligned}
	\end{equation*}
    We denote by $L^{k,Lip}(\mathcal{Z})$ the $\mathcal{C}^{k,Lip}$-norm of $\phi$ for the $\mathcal{C}^{k, Lip}$ domain $\mathcal{Z}$.
    
    The first part of the definition requires that there locally exist a coordinate system such that the boundary coincides locally with the graph of a function whose $k^{th}$ derivative (or the function itself if $k=0$) is Lipschitz. The two inclusions in the second part of the definition guarantee that the interior and the exterior of the domain are nonempty and locally contain a cone around any point of the boundary. They rule out, for example, the domain $\{(x,y) \, : \, x>0,\, |y| < x^2\}\subset \mathbb{R}^2$. Informally, a domain $\mathcal{Z}$ is Lipschitz if its boundary $\partial \mathcal{Z}$ can be viewed locally as the graph of a Lipschitz function for some coordinate system. 

    \item An open bounded subset $\mathcal{Z}$ of $\mathbb{R}^d$ has the {uniform cone property} if for every $x \in \partial \mathcal{Z}$, there exists a neighborhood $V$ of $x$ in $\mathbb{R}^n$ and new coordinates $\{y_1,\dots,y_n\}$ such that (i) $V$ is a hypercube in the new coordinates: $V = \{(y_1,\dots,y_n)\mid -a_i < y_i < a_i,\; 1 \le i \le n\}$, and (ii) $y - z \in \mathcal{Z}$ whenever $y \in \partial \mathcal{Z} \cap V$ and $z \in C$, where $C$ is the open cone $\{z = (z', z_n)\mid (\cot \theta)\,|z'| < z_n < h\}$  for some $\theta \in (0,\pi/2)$) for some $h > 0$. 
    
    We denote by $h(\mathcal{Z})$ and $\theta(\mathcal{Z})$ the parameters of the cone condition. A domain in $\mathbb{R^d}$ has the (uniform) cone property if and only if it is Lipschitz. The parameter of the cone conditions can then be chosen as a function of $L^{0,Lip}(\mathcal{Z})$ only. (See, e.g., Theorem 1.2.2.2 in \cite{grisvard2011elliptic}).

    \item For $\Lambda \geq \lambda \geq 0$, Pucci's extremal operators $\mathcal{P}^{+}_{\lambda, \Lambda}$ and $\mathcal{P}^{-}_{\lambda, \Lambda} : \mathbb{S}^d(\mathbb{R}) \to \mathbb{R}$ are defined by
	\begin{equation*}
	   \mathcal{P}^{+}_{\lambda, \Lambda}(M) = \underset{\lambda I \leq_{\mathbb{S}} B \leq_{\mathbb{S}} \Lambda I}{\sup } \, tr(BM) \, \, \text{ and } \, \, \mathcal{P}^{-}_{\lambda, \Lambda}(M) =  \underset{\lambda I \leq_{\mathbb{S}} B \leq_{\mathbb{S}} \Lambda I}{\inf } \, tr(BM),
	\end{equation*}
    where the order $\geq_{\mathbb{S}}$ is given by $M \geq_{\mathbb{S}} N$ if and only if $M-N$ is positive semi-definite.

    \item We will often use the letter $C$ as a generic letter to denote different bounds appearing in various estimates, and which can be explicitly computed in terms of primitives of the problem. 
\end{itemize}

\subsection{Solution concepts and definitions}\label{subsec:concepts}

\begin{definition}\label{def:strongsolution}
A function $u$ is an $L^p$-subsolution (respectively, $L^p$-supersolution) of~\eqref{eq:obstacle} if the following conditions hold: (i) $u\in {W}^{1,2,p}_{loc}(\mathcal{Y}) \cap \mathcal{C}^0(\bar{\mathcal{Y}})$, (ii) $u \leq b$ (respectively, $\geq b$) on $\partial {\mathcal{Y}}$, and (iii)
	\[\max\left\{ \,  u_t+ F\left(t, x, u, u_x, u_{xx} \right), \,  g - u \right\} \geq 0 \, \,(\text{respectively, } \leq 0) \text{ a.e. on } \mathcal{Y}.\]	
    $u$ is an $L^p$-solution if it is both an $L^p$-subsolution and an $L^p$-supersolution.
\end{definition}
We will also use viscosity solutions in stating and applying the comparison principle of Section~\ref{subsec:comparison}. There are several concepts of viscosity solutions depending on the set of test functions used in the definition. We use  test functions in $W^{1,2,p}\left(\mathcal{Y}\right)$, which corresponds to what is sometimes called ``$L^p$-viscosity solutions.''\footnote{When $F$ is continuous, one can take the test functions to be in $\mathcal{C}^{1,2}$. In this case, the solution is called a $\mathcal{C}$-viscosity solution. In this case, the concepts of $\mathcal{C}$-viscosity solutions and of viscosity solutions as we defined them coincide (Lemma 2.9 in \cite{crandall2000lp}).}
\begin{definition}\label{def:Lpviscositysolution}
A function $u: \bar{\mathcal{Y}} \to \mathbb{R}$ is a viscosity subsolution (respectively, supersolution) of \eqref{eq:obstacle} if (i) it is continuous, (ii) $u \leq b$ (respectively, $\geq b$) on $\partial {\mathcal{Y}}$, and (iii) for all $(t_0, x_0) \in \mathcal{Y}$ and all $\varphi \in {W}^{1,2,p}_{loc}(\mathcal{Y})$ such that $u - \varphi$ (respectively $\varphi-u$) has a maximum at $(t_0, x_0)$, one has
\begin{equation*}
    \begin{aligned}
        \underset{(t,x) \to (t_0, x_0)}{\mathrm{ess}\limsup} &\,  \max\left\{ \, \varphi_t + F\left(t, x, u, \varphi_x, \varphi_{xx}\right), \, g - u \right\}  \geq 0 \\
	\Big(\text{respectively, }  \underset{(t,x) \to (t_0, x_0)}{\mathrm{ess}\liminf} &\,  \max\left\{ \, \varphi_t + F\left(t, x, u, \varphi_x, \varphi_{xx}\right), \, g - u \right\} \leq 0\Big).
    \end{aligned}
\end{equation*}
$u$ is a viscosity solution if it both a viscosity subsolution and supersolution.
\end{definition}

\section{Main Result}\label{sec:mainresult}

Our main result is derived under the following assumptions.
\begin{assumption}\label{assumption:lipschitzdomain}
    $\mathcal{X}$ is $\mathcal{C}^{0, Lip}$.
\end{assumption}
This assumption allows domains $\mathcal{X}$ with non-smooth boundaries. As noted in the Introduction, allowing for such domains is particularly useful in  economic and financial applications, where domains are often equal to a positive orthant (e.g., in the case of price vectors), a unit cube, or a simplex (e.g., when $x$ is a probability distribution representing the decision maker's belief about a state of the world with $d$ possible values). These domains all fail to be smooth, but they all satisfy the Lipschitz property. 

The next two assumptions concern the exogenous boundary payoff function $b$ and the stopping payoff function $g$. 
\begin{assumption}\label{assumption:bcontinuity}
${b}: \bar{\mathcal{Y}} \to \mathbb{R}$ is in $\mathcal{C}^0\left( \bar{\mathcal{Y}} \right)$.\end{assumption}
Assumption \ref{assumption:bcontinuity} implies that $b$ is bounded on $\bar{\mathcal{Y}}$, a property that we will use when computing estimates.

\begin{assumption}\label{assumption:gregularity}
    $g = \underset{a\in A}{\sup} \; g^{a}$ where $A$ is a separable topological space and the functions $g^a$ have the following properties: (i) 
    $g^a \in W^{1,2,p}(\mathcal{Y})$ for all $a$ and $\underset{a \in A}{\sup } \, \left\| g^a \right\|_{W^{1,2,p}\left(\mathcal{Y}\right)} <\infty$, (ii) $g^a \leq b$ on $\partial \mathcal{Y}$,  
    and (iii) the map $a\mapsto g^a$ is continuous from $A$ to $C^{0}\left(\bar{\mathcal{Y}}\right)$.\end{assumption}

Part (iii) of Assumption~3 guarantees that the obstacle can be approximated by the maximum of finitely many $W^{1,2,p}$-functions, which will be key for the proof of Theorem \ref{theorem:existenceobstacleproblemirregulardomain}.\footnote{Part (iii) is sufficient but need not be necessary. The proof requires only the existence of an approximating sequence $(g^n = \underset{a\in A_n}{\max} \, g^{a})_n$, with $A_n$ finite for all $n$.} Moreover, it ensures that $g$ is measurable as it can be written as the supremum of a family of measurable function over a countable dense subset of $A$.

The next assumption is a \textbf{structure condition}.
\begin{assumption}\label{assumption:structurecondition}
There exist $\lambda, \Lambda >0$ and moduli of continuity $\omega_1$ and $\omega_2$ (i.e., nondecreasing continuous functions on $\mathbb{R}_+$ with $\omega_j(0) = 0$, $j=1,2$) such that, for almost all $(t,x) \in \mathcal{Y}$, we have
    \begin{equation}\label{eq:SC}
        \begin{cases}
            & \mathcal{P}_{\lambda, \Lambda}^-\left(M-\tilde{M}\right) - \omega_{1}\left( \right|r -\tilde{r} \left|\right) - \omega_{2} \left(\right|q -\tilde{q} \left| \right) \\
            & \quad \quad \leq  F\left(t, x, r, q, M\right) - F\left(t, x, \tilde{r}, \tilde{q}, \tilde{M} \right)  \\
            &\quad \quad \quad \quad \leq \mathcal{P}_{\lambda, \Lambda}^+\left(M-\tilde{M}\right) + \omega_{1}\left( \right|r -\tilde{r} \left|\right) +\omega_{2} \left(\right|q-\tilde{q} \left| \right).
    	\end{cases}
     \end{equation}
     for all $M,\tilde{M} \in \mathbb{S}^d$, $r, \tilde{r} \in \mathbb{R}$, and $q, \tilde{q} \in \mathbb{R}^{d}$, where $\mathcal{P}^-_{\lambda, \Lambda}$ and $\mathcal{P}^+_{\lambda, \Lambda}$ are the Pucci extremal operators. Moreover, we assume that there exists $R>0$ such that $\omega_1(r) +\omega_2(q) \leq R(1 + r + q)$ for all $r,q \in \mathbb{R}_+$.
\end{assumption}
Assumption~\ref{assumption:structurecondition} implies that $F$ is uniformly elliptic -- i.e., that $F$ is (strongly) increasing in $M$. Moreover, it  that $F(t,x, r,p,M)$ is continuous in $r$, $p$, and $M$ uniformly over $(t,x)$, and that $F$ grows at most linearly in the value and first derivatives of the value function. 

In the literature, stronger regularity conditions on $F$ are typically imposed to obtain existence and regularity of solutions. In particular, $F$ is generally assumed to be Lipschitz continuous in $(r,q,M)$ and monotone in $r$.  See, e.g., \cite{crandall2000lp} or \cite{dong2013fully} By contrast, Assumption \ref{assumption:structurecondition} only requires continuity in $(r,q,M)$ and does not impose any monotonicity condition with respect to $r$,\footnote{To establish the uniqueness part of our result, and only this part, we do impose a weak monotonicity condition (Assumption 8).} hence accommodating a larger class of operators.

\begin{assumption}\label{assumption:convexity}$F(t,x, 0,0,M)$ is convex in $M$ for all $(t,x)$.\footnote{When $F$ is concave case, our results remain valid (upon modifying Assumption \ref{assumption:strictmonotonicity} appropriately), since $u$ solves $\max\{ g-v, F(t,x,v,v_x,v_{xx}) + v_t\} = 0$ if and only if it solves $\min\{ v - g , - F(t,x,v,v_x,v_{xx}) - v_t\} = 0$.}
\end{assumption}

The following Vanishing Mean Oscillation (VMO) assumption is also imposed. Define the VMO modulus of $F$ by
\begin{equation*}
    \eta_F(\delta) = \underset{Q \subset \mathcal{Y} \ : \ \left|Q \right| \leq \delta }{\sup} \, \frac{1}{\left|Q\right|} \int_{Q} \underset{M \in \mathbb{S}^d \setminus \{0\}}{\sup }\, \frac{\left| F\left(t, x, 0, 0, M\right) - F_{Q}(M)\right|}{\left\| M \right\|} dt dx,
\end{equation*}
where
\begin{equation*}
    F_{Q}(M) = \frac{1}{\left|Q\right|} \int_{Q} F\left(t, x, 0, 0, M\right) dt dx.
\end{equation*}

\begin{assumption}\label{assumption:VMO}
    $\frac{F(t,x, 0,0, M)}{\left\|M\right\|}$ is in the Sarason class VMO (see \cite{sarason1975functions}), i.e., $\underset{\delta\to 0^+}{\lim \ } \eta_F(\delta) = 0$.
\end{assumption}
Assumption 6 imposes a vanishing mean oscillation (VMO) condition on the dependence of $F$ on the second-order term. Intuitively, it requires that, at small scales, the operator $F(t,x,0,0, M)$ behaves alike one that does not depend on $(t,x)$. In particular, it is satisfied if $\frac{F(t,x,0,0, M)}{\left\|M\right\|}$ is continuous (uniformly in $M$). If $F$ were linear -- i.e., $F(t,x,0,0,M) = {tr}(\sigma\sigma'(t,x) M)$ -- then Assumption 6 reduces to requiring that the matrix-valued coefficient $\sigma\sigma'(t,x)$ belongs to VMO -- a standard condition in the regularity theory of parabolic equations.

Assumption \ref{assumption:VMO} is needed to ensure that the weak derivatives of solutions to canonical equations related to obstacle problem exist a.e. and belong to $L^p$, by guaranteeing that the coefficient of the highest-order derivative in the linearized version of the operator $F$ is regular enough in $(t,x)$. As pointed out in \cite{dong2020recent}, Assumption \ref{assumption:VMO} is the weakest known assumption (even in the linear case) under which a solution in $W^{1,2,p}$ always exists for Dirichlet problems. Moreover, counter-examples in \cite{meyers1963p} or \cite{ural1967невозможности} suggest that this condition is difficult to relax.

The next assumption is standard and necessary to bound the  $W^{1,2,p}$-norm of solutions to \eqref{eq:obstacle}.
\begin{assumption}\label{assumption:growthconditions}
    There exists ${G} \in {L}^{p}\left(\mathcal{Y}\right)$ such that, for all $\left(t,x \right) \in \mathcal{Y}$,
    \begin{equation*}
        \left| F\left( t,x, 0, 0, 0 \right) \right|\leq  G(t,x).
    \end{equation*}
\end{assumption}
The final assumption is used to apply a comparison principle to \eqref{eq:obstacle}. It is not needed for our existence results, but it is essential to guarantee the uniqueness of a solution to \eqref{eq:obstacle}.
\begin{assumption}\label{assumption:strictmonotonicity}
    There exists $\kappa>0$ such that, for all $(t,x) \in \mathcal{Y}$, $q \in \mathbb{R}^{d}$, and $M \in \mathbb{S}^d$,
    \begin{equation*}
        r \to F(t,x,r,q,M) - \kappa r
    \end{equation*}
    is strictly decreasing.
\end{assumption}

Our main result establishes the existence, uniqueness, and regularity of solutions to\footnote{The structure condition \eqref{eq:SC} guarantees that any solution of $u_t + F[u]=0$ with $u =b$ on $\partial \mathcal{Y}$ is bounded below. Therefore, by choosing $g$ small enough, Theorem~1 also guarantees the existence, uniqueness, and regularity of the solution without obstacles.}  \eqref{eq:obstacle}.

\begin{theorem}\label{theorem:existenceobstacleproblemirregulardomain}
    Suppose that Assumptions \ref{assumption:lipschitzdomain}--\ref{assumption:growthconditions} hold for some\footnote{The parameter $p$ appears in Assumptions~3 and 7.} $p \in (d+2, \infty)$. Then \eqref{eq:obstacle} has an $L^p$-solution $u$. Moreover, for all compact $\mathcal{Y}' \subset \mathcal{Y}$, there exists $C = C\left( d, p, \lambda, \Lambda, R, \eta_F, diam\left(\mathcal{X}\right), T, dist(\mathcal{Y}', \mathcal{Y}) \right) \in \mathbb{R}_+$ such that
\begin{equation}\label{eq:W12pestimatesirregulardomain}
    \left\| u \right\|_{W^{1,2,p}\left(\mathcal{Y'}\right)} \leq C \left( 1 + \underset{a\in A}{\sup}\, \left\| g^a \right\|_{W^{1,2,p} \left( \mathcal{Y}\right)} + \left\| G \right\|_{L^{p}\left( \mathcal{Y}\right)} + \left\|b \right\|_{L^{\infty}\left( \mathcal{Y} \right)} \right).
\end{equation}

If, in addition, Assumption \ref{assumption:strictmonotonicity} holds, $u$ is the unique $L^p$-solution of \eqref{eq:obstacle}.
\end{theorem}

\begin{remark} The $W^{1,2,p}$-estimate~\eqref{eq:W12pestimatesirregulardomain} is independent of the moduli of continuity $\omega_1$ and $\omega_2$ introduced in Assumption~\ref{assumption:structurecondition}.\end{remark}

\begin{corollary}\label{corollary:smoothness}
     Suppose that Assumptions \ref{assumption:lipschitzdomain}--\ref{assumption:strictmonotonicity} hold for some $p \in (d+2, \infty)$. Then, the unique $L^p$-solution $u$ of \eqref{eq:obstacle} is continuously differentiable with respect to the space variable $x$ on $\mathcal{Y}$.
\end{corollary}

For stopping problems, the solution of \eqref{eq:obstacle} coincides with the value function. Corollary \ref{corollary:smoothness} then (i) implies that smooth pasting in space holds, and (ii) confirms that the decision maker never stops at a kink of $g$, as shown by \cite{decamps2006irreversible} in the specific problem of irreversible investment.

\vspace{0.3em}
\proof{Proof.}
    By Theorem \ref{theorem:existenceobstacleproblemirregulardomain}, the unique solution is in $W^{1,2,p}_{loc}(\mathcal{Y})$. By a Morrey-Sobolev embedding type theorem (Lemma 3.3, page 80 in \cite{ladyzenskaja1968linear}), $W^{1,2,p}_{loc}\left(\mathcal{Y}\right)$ is embedded in $\mathcal{C}^{0,1}\left( \mathcal{Y} \right)$. 
    
    \hfill $\square$
\endproof

\section{Application: Robustly Investing in Alternative Projects}\label{sec:application}

In this section, we study an irreversible investment problem with model uncertainty and multiple investment options, inspired by \cite{hansen2001robust} and \cite{decamps2006irreversible}, to illustrate our results. 

A risk-neutral decision maker can invest in one of $I$ projects, $a \in \{1, \dots, I\}$, whose return upon investment depends on a vector of primitives $X$. The decision-maker must choose (i) when, before a deadline $T < \infty$, to invest, and (ii) in which project $a\in \{1,\dots, I\}$ to invest. The state variable $X$ takes values in $\mathcal{X}\subset \mathbb{R}^d$ and captures the evolution of the characteristics of the projects, e.g., the dynamics of output prices. It evolves according to the diffusion equation:
\begin{equation}\label{eq:sdebenchmark}
    X_t = X_0 + \int_0^t \mu(s, X_s) ds + \int_0^t \sigma(s, X_s) dB_s,
\end{equation}
where $B$ is an $n$-dimensional standard Brownian motion. In order to apply our main Theorem \ref{theorem:existenceobstacleproblemirregulardomain}, we assume that the set $\mathcal{X}$ is Lipschitz and bounded, and hence, so is the parabolic domain $\mathcal{Y} = [0,T) \times \mathcal{X}$. Therefore Assumption \ref{assumption:lipschitzdomain} hold. Moreover, to ensure that Assumptions \ref{assumption:structurecondition}, \ref{assumption:VMO}, and \ref{assumption:growthconditions} are verified, we assume that the drift $\mu: \mathcal{Y} \to \mathbb{R}^d$ is measurable and bounded, and that the volatility $\sigma:\mathcal{Y}\to \mathbb{R}^{d\times n}$ is a bounded VMO function and is uniformly nondegenerate: there exist $\lambda, \Lambda > 0$ such that
\begin{equation*}
    \lambda I_d \leq_{\mathbb{S}} \sigma\sigma'(t,x) \leq_{\mathbb{S}} \Lambda I_d
\quad \text{ for all } (t,x) \in \mathcal{Y}.
\end{equation*}
Finally, each project $a$ yields a payoff $g^a(t,x)$ when chosen at time $t$ in state $x$, where $g^a \in W^{1,2,p}(\mathcal{Y})$. So, the decision maker stopping payoff is $\underset{a \in \{1, \dots, I\}}{\max} \ g^a$, which guarantees that Assumption \ref{assumption:gregularity} holds; and we assume that the DM discounts future payoffs at rate $r>0$ and must make a decision at time $\tau_{\mathcal{Y}} = \inf\left\{s\geq t \ : \ (s,X_s) \not \in \mathcal{Y}\right\}$ when $(t,x)$ first exits the domain $\mathcal{Y} = [0,T) \times \mathcal{X}$. So, her boundary payoff is also equal to $\underset{a \in \{1, \dots, I\}}{\max} \ g^a$, which is continuous. Hence, Assumption \ref{assumption:bcontinuity} is verified.

Following \cite{hansen2001robust}, we consider the robust investment problem. The decision maker is concerned about model misspecification. That is, the decision maker considers the possibility that the dynamics of $X$ could be different from the baseline model \eqref{eq:sdebenchmark}. Specifically, she treats the benchmark model \eqref{eq:sdebenchmark} as an approximation, and considers a family of alternative models that are statistically difficult to distinguish from it:
\begin{equation}\label{eq:sdedistorted}
    X_t = X_0 + \int_0^t \left( \mu(s, X_s) + \sigma(s, X_s) h_s \right) ds + \int_0^t \sigma(s, X_s)\, dB_s,
\end{equation}
where $h$ is progressively measurable taking values in $\mathbb{R}^n$ that satisfies $|h_s| \leq H$ for some constant $H \geq 0$. Because she does not know which model represents the true evolution of $X$, the decision-maker optimizes over the worst case scenario for the dynamics. However, to account for the fact that larger model departures are less likely, such departures are ``penalized'' by an in increase in the objective function (typically assumed proportional to the relative entropy between the probability distribution of the benchmark model \eqref{eq:sdebenchmark} and the probability distribution of the distorted models \eqref{eq:sdedistorted}), which makes them less relevant for the worst case scenario, other things equal.

\cite{hansen2001robust} observed that solving the robust control problem is therefore equivalent to solving a zero-sum game between the decision-maker and an adversarial nature. In particular, in the robust investment problem we examine, the decision-maker behavior is identical to her behavior when she plays a zero-sum game against an adversarial nature who can choose any weak control $h$ -- i.e., controls $h$ such that \eqref{eq:sdedistorted} admits a weak solution\footnote{Any Markovian control satisfies this condition. Moreover, when $\sigma$ is VMO and uniformly nondegenerate, as we assume, the Zvonkin transformation method together with the results of this paper guarantees that \eqref{eq:sdedistorted} admits a strong solution for any Markovian control $h(t,x)$. A formal proof is beyond the scope of this paper.} -- and is penalized through a quadratic cost $\frac{\kappa}{2} |h_s|^2$. Thus, the decision-maker chooses a stopping time to maximize her worst-case value: For all $(t,x) \in \mathcal{Y}$,
\begin{equation*}
    V(t,x) = \sup_{\tau} \inf_{h:\, |h_s|\leq H} \mathbb{E}_{(t,x)}\bigg[
    \int_t^{\tau \wedge \tau_{\mathcal{Y}}} e^{-r(s-t)} \frac{\kappa}{2} |h_s|^2 \, ds 
    + e^{- r (\tau \wedge \tau_{\mathcal{Y}}-t)} 
    g(\tau \wedge \tau_{\mathcal{Y}}, X_{\tau \wedge \tau_{\mathcal{Y}}})
    \bigg],
\end{equation*}
where $\mathbb{E}_{(t,x)}$ denotes expectation conditional on $X_t = x$ under the dynamics \eqref{eq:sdedistorted}. The parameters $\kappa$ and $H$ measure ambiguity aversion. When $\kappa = \infty$ or $H = 0$, the problem reduces to the standard irreversible investment problem, as in \cite{dixit1993choosing} or \cite{decamps2006irreversible}. On the other hand, the decision maker hedges against increasingly large misspecifications as $H \to \infty$ and $\kappa \to 0$.

To solve the decision maker’s robust investment problem, we compute the (upper) value of the game, which we expect to solve the Hamilton–Jacobi–Bellman–Isaacs equation:
\begin{equation}\label{eq:HJBI}
    \begin{cases}
        \max\bigg\{ u_t + \underset{h \in \mathbb{R}^n : \left| h \right| \leq H}{\inf} \ \left[ \frac{\kappa}{2} |h|^2 
        + (\mu(t,x) + \sigma(t,x) h ) u_x \right] 
        + \frac{1}{2} tr\left( \sigma\sigma'(t,x) u_{xx} \right) - r u, \\
        \qquad\qquad\ \underset{a\in \{1,\dots, I\}}{\max} g^a - u \bigg\} = 0 
        \quad \text{in } \mathcal{Y}, \\
        u = \underset{a\in \{1,\dots, I\}}{\max} g^a \quad \text{on } \partial \mathcal{Y},
    \end{cases}
\end{equation}
provided that a sufficiently regular solution exists (see, e.g., \cite{fleming1989existence}).

However, establishing the existence of such a sufficiently regular solution (so a verification argument applies) can be challenging. Existing results do not apply as the obstacle \(\max_{a \in \{1,\dots,I\}} g^a\) generally fails to belong to \(W^{1,2,p}(\mathcal{Y})\) and the Hamiltonian operator is only quasi-linear. By contrast, the tools developed in this paper are designed precisely to handle such irregular obstacles. This allows us to show that the value of the robust investment problem coincides with the solution to \eqref{eq:HJBI}.
\begin{proposition}\label{prop:robust_stopping}
    Let $p>d+2$. The value $V(t,x)$ is the \textbf{unique $L^p$-solution} of \eqref{eq:HJBI}. In particular, $V$ is continuously differentiable with respect to the space variable $x$.

    Finally, the stopping time $\tau^{\star} = \inf\left\{s \geq t \: \  u(s, X_s) = \underset{a\in \{1,\dots, I\}}{\max} g^a(s, X_s) \right\}$ and decision $a^{\star} \in \arg \underset{a\in \{1,\dots, I\}}{\max} g^a(\tau^{\star}, X_{\tau^{\star}})$ are robustly optimal in the continuation problem $V(t,x)$. 
\end{proposition}
\proof{Proof.} The proof proceeds in two steps. First, we apply Theorem \ref{theorem:existenceobstacleproblemirregulardomain} to establish the existence of an $L^p$-solution to \eqref{eq:HJBI}. Second, we show that any $L^p$-solution of \eqref{eq:HJBI} coincides with the value function $V$ by a standard verification argument based on Krylov–It\^{o}'s lemma.

\textbf{Step 1: Existence.} We first show that a $L^p$-solution of \eqref{eq:HJBI} exists for all $p>d+2$. 
    
    This is a consequence of Theorem \ref{theorem:existenceobstacleproblemirregulardomain}, which applies to \eqref{eq:HJBI}. To see this, let $p>d+2$. Note that Assumptions \ref{assumption:lipschitzdomain}, \ref{assumption:bcontinuity}, and \ref{assumption:gregularity} are satisfied by definition of the domain $\mathcal{X}$ and stopping payoff (which we also take to be the boundary value). Next, we show that $F(t,x, u, u_x, u_{xx}) =  \underset{h : \left| h \right| \leq H}{\inf} \ \frac{\kappa}{2} (h \cdot h ) -r u + \left( \mu(t, x) + h \sigma(t, x) \right) u_x + \frac{1}{2} tr\left( \sigma (t,x) \sigma'(t,x) u_{xx} \right)$ satisfies Assumptions \ref{assumption:structurecondition}, \ref{assumption:convexity}, \ref{assumption:VMO}, and \ref{assumption:growthconditions}. First, observe that
    \begin{equation*}
        \begin{aligned}
         & F \left(t, x, u, u_x, u_{xx}\right) - F\left(t, x, \tilde{u}, \tilde{u}_x, \tilde{u}_{xx} \right) \\
         & = \underset{h \in \mathbb{R}^n : \left| h \right| \leq H}{\inf} \ \frac{\kappa}{2} (h \cdot h ) -r u + \left( \mu(t, x) +  \sigma(t, x) h \right) u_x + \frac{1}{2} tr\left( \sigma \sigma'(t,x) u_{xx} \right) \\
         & \qquad - \underset{h \in \mathbb{R}^n : \left| h \right| \leq H}{\inf} \ \frac{\kappa}{2} (h \cdot h ) -r \tilde{u} + \left( \mu(t, x) +  \sigma(t, x) h \right) \tilde{u}_x + \frac{1}{2} tr\left( \sigma \sigma'(t,x) \tilde{u}_{xx} \right) \\
         & \leq \underset{h \in \mathbb{R}^n : \left| h \right| \leq H}{\sup} \ \frac{\kappa}{2} (h \cdot h ) -r (u - \tilde{u}) + \left( \mu(t, x) +  \sigma(t, x) h \right) (u_x - \tilde{u}_x) + \frac{1}{2} tr\left( \sigma \sigma'(t,x) (u_{xx} - \tilde{u}_{xx}) \right) \\
         & \leq \mathcal{P}_{\lambda, \Lambda}^+\left(u_{xx}-\tilde{u}_{xx}\right) + r \left|u -\tilde{u} \right| + \left( \left\|\mu \right\|_{L^{\infty}(\mathcal{Y})} + H \left\|\sigma \right\|_{L^{\infty}(\mathcal{Y})} \right)\left|q-\tilde{q} \right| .
         \end{aligned}
    \end{equation*}
    which proves the right inequality in \eqref{eq:SC}. The left inequality is established similarly. So, Assumption \eqref{assumption:structurecondition} holds. Second, since the adversarial control does not affect the volatility, $F$ is linear in $u_{xx}$, so Assumption \ref{assumption:convexity} holds. Moreover, $\left| F(t,x, 0, 0, M) - F(t',x', 0, 0, M)\right| \leq \frac{1}{2} \left| tr\left( \left( \sigma (t,x) \sigma'(t,x)  - \sigma (t',x') \sigma'(t',x') \right) M \right)\right| $, which satisfy Assumption \ref{assumption:VMO}, since the volatility function $\sigma$ is a bounded VMO function and the product of VMO functions is VMO. Finally, $\left| F(t,x, 0, 0, 0) \right| \leq \frac{\kappa}{2} H^2$, and, hence, $F$ satisfies Assumption \ref{assumption:growthconditions}. 
    
    Therefore, Theorem \ref{theorem:existenceobstacleproblemirregulardomain} applies and \eqref{eq:HJBI} has a $L^p$-solution.

    \textbf{Step 2: Verification.} We show that any $L^p$-solution of \eqref{eq:HJBI} coincides with the value function. 
    
    Let $u$ be a $L^p$-solution of \eqref{eq:HJBI}. Let $(t,x) \in \bar{\mathcal{Y}}$. If $(t,x) \in \partial \mathcal{Y}$, then $u = \underset{a}{\max} \ g^a = V$. If, instead, $(t,x) \in \mathcal{Y}$, consider the adversarial control $h^{\star}_s = h^{\star}(s,X_s)$, where $(t,x) \to h^{\star}(t,x)$ is a measurable selection that solves the inner minimization problem in\footnote{Such a selection exists by, e.g., the Measurable Maximum Theorem 18.19 in \cite{aliprantis2006infinite}.} \eqref{eq:HJBI} and the stopping stopping time $t\leq \tau \leq \tau_{\mathcal{Y}_n} = \inf\left\{ s \geq t \ : \ (s, X_s) \not \in [0,T) \times \mathcal{X}_n\right\}$, where $\mathcal{X}_n$ is an increasing sequence of smooth domains such that $\bigcup_n \mathcal{X}_n = \mathcal{X}$, as defined in the proof of Theorem \ref{theorem:existenceobstacleproblemirregulardomain}. Since $u \in W^{1,2,p}_{loc}(\mathcal{Y}) \cap \mathcal{C}^0(\bar{\mathcal{Y}})$, It\^{o}-Krylov's formula (Theorem 2.10.1 in \cite{krylov2008controlled}) implies that
    \begin{equation*}
        \begin{aligned}
        u(t,x) & = \mathbb{E}_{(t,x)} \bigg[ e^{-r (\tau-t)} u(\tau, X_{\tau}) - \int_{t}^{\tau}e^{- r(s-t)} \bigg( u_t(s, X_s) +  \left( \mu(s, X_s) + \sigma(s, X_s) h^{\star}_s \right) u_x(s, X_s) \\
        & \qquad \qquad + \frac{1}{2} tr\left( \sigma (t,x) \sigma'(t,x) u_{xx}(s, X_s)\right)  - ru(s, X_s)  \bigg) ds\bigg].
        \end{aligned}
    \end{equation*}
    Since $u$ is an $L^p$-solution of \eqref{eq:HJBI}, by definition of $h^{\star}_s$, for almost all $(t',x') \in \mathcal{Y}$,
    \begin{equation*}
        \begin{aligned}
        u_t(t',x') & +  \left( \mu(t', x') + \sigma(t', x') h^{\star}(t',x') \right) u_x(t', x') \\
        &  \qquad + \frac{1}{2} tr\left( \sigma (t',x') \sigma'(t',x') u_{xx}(t', x')\right)  - ru(t', x') \leq - \frac{\kappa}{2} \left(h^{\star}(t',x') \cdot h^{\star}(t',x')\right).
        \end{aligned}
    \end{equation*}
    Therefore, 
    \begin{align}\label{eq:proofverification3}
        u(t,x) \geq \mathbb{E}_{(t,x)} & \bigg[  e^{-r (\tau-t)} u(\tau, X_{\tau}) + \int_{t}^{\tau} e^{-r(s-t)} \frac{\kappa}{2} \left(h^{\star}(s,X_s) \cdot h^{\star}(s,X_s)\right)ds \bigg].
    \end{align}
    Moreover, since $u$ is a $L^p$-solution of \eqref{eq:HJBI}, $u \geq \underset{a\in \{1,\dots, I\}}{\max} g^a$ on $\bar{\mathcal{Y}}$. Thus,
    \begin{equation*}
        \begin{aligned}
        u(t,x) & \geq \mathbb{E}_{(t,x)} \bigg[  e^{-r (\tau-t)} \underset{a\in \{1,\dots, I\}}{\max} g^a(\tau, X_{\tau}) + \int_{t}^{\tau} e^{-r(s-t)} \frac{\kappa}{2} \left(h^{\star}(s,X_s) \cdot h^{\star}(s,X_s)\right)ds \bigg] \notag \\
        & \geq \underset{h : \left|h_s \right| \leq H }{\inf } \ \mathbb{E}_{(t,x)} \bigg[  e^{-r (\tau-t)} \underset{a\in \{1,\dots, I\}}{\max} g^a(\tau, X_{\tau}) + \int_{t}^{\tau} e^{-r(s-t)} \frac{\kappa}{2} \left(h_s) \cdot h_s\right)ds \bigg].
        \end{aligned}
    \end{equation*}
    The right-hand side of the above equation is the expected worst-case payoff for the arbitrary stopping time $t\leq \tau \leq \tau_{\mathcal{Y}}$. Letting $n\to \infty$, and taking the supremum then shows that $u(t,x) \geq V(t,x)$. Since $(t,x) \in \mathcal{Y}$ were arbitrary, $u \geq V$ on $\bar{\mathcal{Y}}$.

    Finally, when $\tau = \tau^{\star} \wedge \tau_{\mathcal{Y}_n}$ and $h = h^{\star}$, the inequality \eqref{eq:proofverification3} holds with equality, i.e., 
    \begin{equation*}
        \begin{aligned}
        u(t,x) & = \mathbb{E}_{(t,x)} \bigg[  e^{-r (\tau^{\star}\wedge \tau_{\mathcal{Y}_n}-t)} u(\tau^{\star} \wedge \tau_{\mathcal{Y}_n}, X_{\tau^{\star} \wedge \tau_{\mathcal{Y}_n}}) \\
        & \qquad \qquad + \int_{t}^{\tau^{\star} \wedge \tau_{\mathcal{Y}_n}} e^{-r(s-t)} \frac{\kappa}{2} \left(h^{\star}(s,X_s) \cdot h^{\star}(s,X_s)\right)ds \bigg].
        \end{aligned}
    \end{equation*}
    Hence, the dominated convergence theorem (using that $u \in \mathcal{C}^0(\mathcal{\bar{Y}})$) and the definition of $\tau^{\star}$ and $a^{\star}$ yield
    \begin{equation*}
        \begin{aligned}
        u(t,x) & = \mathbb{E}_{(t,x)} \bigg[  e^{-r (\tau^{\star}-t)} u(\tau^{\star}, X_{\tau^{\star}}) + \int_{t}^{\tau^{\star}} e^{-r(s-t)} \frac{\kappa}{2} \left(h^{\star}(s,X_s) \cdot h^{\star}(s,X_s)\right)ds \bigg] \\
        & = \mathbb{E}_{(t,x)} \bigg[  e^{-r (\tau^{\star}-t)} \underset{a\in \{1,\dots, I\}}{\max} g^a(\tau^{\star}, X_{\tau^{\star}}) + \int_{t}^{\tau^{\star}} e^{-r(s-t)} \frac{\kappa}{2} \left(h^{\star}(s,X_s) \cdot h^{\star}(s,X_s)\right)ds \bigg]
        \end{aligned}
    \end{equation*}
    So, $\tau^{\star}$ and $a^{\star}$ are robustly optimal and $V(t,x) = u(t,x)$.
    \hfill $\square$

Proposition \ref{prop:robust_stopping} shows that the value function is continuously differentiable in space despite the possible kinks in the stopping payoff \(\max_{a} g^a\). Therefore, the optimal stopping region cannot include any non-smooth indifference curve. In economic terms, the decision maker never invests at a point where they are indifferent between projects and where waiting an infinitesimal amount of time would break this indifference. This confirms the insight of \cite{decamps2006irreversible}, who show, using probabilistic methods, that the optimal stopping region can be disconnected in a (non-robust) irreversible investment problem due to the non-differentiability of the payoff. 

Finally, Proposition \ref{prop:robust_stopping} also provides insights into the behavior of a decision maker concerned about model misspecification. For example, we obtain the following corollary using the comparison principle (Proposition \ref{prop:comparison}) derived below.
\begin{corollary}\label{corollary:compstats}
    The stopping region $\mathcal{S} = \left\{ (t,x) \in \bar{\mathcal{Y}} \ : \ u(t,x) = \underset{a\in \{1,\dots, I\}}{\max} g^a(t,x)  \right\}$ expands as the decision maker becomes more ambiguity averse -- i.e., as $\kappa$ decreases and $H$ increases.
\end{corollary}

\section{Proof of Theorem~\ref{theorem:existenceobstacleproblemirregulardomain}}\label{sec:proof}

In this section, we prove our main result, Theorem \ref{theorem:existenceobstacleproblemirregulardomain}, in several steps. We first prove two key auxiliary results in Section \ref{sec:comparisonandstability}: (i) a comparison principle for fully nonlinear parabolic obstacle problems with measurable ingredients (Proposition \ref{prop:comparison}), and (ii) a stability theorem for solutions in $W^{1,2,p}(\mathcal{Y})$ (Theorem \ref{theorem:continuityequation}). Next, in Section \ref{sec:finiteA}, we prove a more restrictive version of Theorem \ref{theorem:existenceobstacleproblemirregulardomain} when the obstacle is the maximum of finitely many $W^{1,2,p}$ functions and both the domain and boundary data are regular. Finally, Section \ref{sec:proofirregulardomains} contains the proof of Theorem \ref{theorem:existenceobstacleproblemirregulardomain}.

\subsection{Comparison and Stability for fully nonlinear parabolic obstacle problems}\label{sec:comparisonandstability}

Section \ref{sec:comparisonandstability} contains two results needed for the proof of Theorem \ref{theorem:existenceobstacleproblemirregulardomain}: (i) A comparison principle that applies to fully nonlinear parabolic obstacle problems with measurable ingredients (Proposition \ref{prop:comparison}), and (ii) a stability theorem for solutions to obstacle problems in $W^{1,2,p}\left(\mathcal{Y}\right)$ (Theorem \ref{theorem:continuityequation}).

We start with a useful well-known observation: Viscosity solutions are weaker than $L^p$-solutions.
\begin{lemma}\label{lemma:relationdifferentsolutionconcepts}
Let $u$ be an $L^p$-subsolution (respectively, supersolution) of \eqref{eq:obstacle}. Then $u$ is a viscosity subsolution (respectively, supersolution) of \eqref{eq:obstacle}.\end{lemma}
Lemma \ref{lemma:relationdifferentsolutionconcepts} has been proved in, e.g., \cite{crandall2000lp} (Proposition 2.11) in a setting similar to ours. It follows easily from the definition of $L^p$-viscosity solution, the monotonicity of the operator $F$ in the Hessian, and Bony's maximum principle (see Theorem 4.2 in \cite{yazhe1985aleksandrov}).

\subsubsection{Comparison Principle}\label{subsec:comparison}
Comparison principles for $L^p$-solutions and viscosity solutions of the Dirichlet problem are well-known in the literature. We provide here a more general result that encompasses obstacle problems and includes operators that are less restrained than those considered in the literature.\footnote{For instance, Proposition 2.10 in \cite{crandall2000lp} gives a comparison principle when $F$ is Lipschitz continuous in $u_x$ and nonincreasing in $u$.} Our comparison principle is valid for a large class of nonlinear parabolic obstacle problems as long as $F$ satisfies the structure condition \eqref{eq:SC} and does not increase too fast with $u$ (as required by Assumption 8).
\begin{proposition}\label{prop:comparison}
    Suppose that $F$ satisfies Assumptions \ref{assumption:structurecondition} and \ref{assumption:strictmonotonicity}.
\begin{itemize}
    \item Let $u$ be an $L^p$-subsolution\footnote{Note that we do not impose here conditions on the functions at the boundary $\partial \mathcal{Y}$; that is, we dispense with requirement (ii) in Definitions 1 and 2 and  allow $u$ and $v$ to take a priori any values on $\partial \mathcal{Y}$. This is because the comparison principle only depends on the condition that $v\geq u$ on $\partial \mathcal{Y}$, which can be used to satisfy Definitions 1 and 2, including requirement (ii), for functions that depend on $u$ and $v$ only through the difference $v-u$.} and $v$ be a viscosity supersolution of
        \begin{equation}\label{eq:comparison1}
            \max\left\{ \,  u_t+ F\left(t, x, u, u_x, u_{xx} \right), \,  g - u \right\} = 0.
        \end{equation}
        If $v \geq u$ on $\partial \mathcal{Y}$, then $v \geq u$ on $\bar{\mathcal{Y}}$.
        \item Let $u$ be an $L^p$-supersolution\footnote{The same observation regarding the absence of a constraint on $\partial \mathcal{Y}$ applies here as in the previous endnote.} and $v$ be a viscosity subsolution of
        \begin{equation}\label{eq:comparison2}
            \max\left\{ \,  u_t+ F\left(t, x, u, u_x, u_{xx} \right), \, g - u \right\} =0.
        \end{equation}
        If $v \leq u$ on $\partial \mathcal{Y}$, then $v \leq u$ on $\bar{\mathcal{Y}}$.
    \end{itemize}
\end{proposition}

The proof of Proposition \ref{prop:comparison} builds on a comparison principle for viscosity solutions of fully nonlinear parabolic Dirichlet problems with continuous operators due to \cite{giga1991comparison}.

\vspace{0.3em}
\proof{Proof.}
Let $u$ be an $L^p$-subsolution and $v$ be a viscosity supersolution of 
\begin{equation*}
        \max\left\{ \,  u_t+ F\left(t, x, u, u_x, u_{xx} \right), \,  g - u \right\} = 0
\end{equation*}
such that $v \geq u$ on $\partial \mathcal{Y}$. We will prove that $v \geq u$ on $\bar{\mathcal{Y}}$. 
    
Let $w = u-v$ and $\mathcal{E} = \left\{(t,x) \in \bar{\mathcal{Y}} \, : \, w > 0\right\}$. We need to show that $\mathcal{E}$ is empty. Suppose by way of contradiction that it is not. We note that $\mathcal{E} \subset \mathcal{Y}$ since $w \leq 0$ on $\partial \mathcal{Y}$. Moreover, $w$ is continuous since both $u$ and $v$ are by construction continuous on $\bar{\mathcal{Y}}$ . Therefore, $\mathcal{E}$ is an open bounded domain. 

We begin by showing that $w$ is a viscosity subsolution of the problem
\begin{equation}\label{eq:proofcomparison1}
\begin{cases}
\max \left\{ \, w_t + H\left(t,x,w, w_x,  w_{xx}\right), \, -w \right\} = 0 \text{ on } \mathcal{E} \\
w = 0 \text{ on } \partial \mathcal{E}
\end{cases}
\end{equation}
where
\begin{equation*}
    	H(t,x,w, w_x,w_{xx}) =  F\left(t, x,  u, u_x, u_{xx}\right) - F\left(t, x, u-w, u_x - w_x, u_{xx} -w_{xx}\right).
\end{equation*}
To see this, we first note that $w = 0$ on $\partial \mathcal{E}$, so condition (ii) in the definition of viscosity subsolutions (which, here, amounts to $w\leq 0$) is satisfied. Next, let $\varphi \in W^{1,2,p}\left(\mathcal{E}\right)$ be such that $w- \varphi$ has a maximum at $(t_0,x_0)$ in $\mathcal{E}$, i.e., $u-\varphi -v$ has a maximum at $(t_0,x_0)$. By construction, $u-\varphi \in W^{1,2,p}_{loc}\left(\mathcal{E}\right)$ and $v$ is a viscosity supersolution of \eqref{eq:comparison1}. Therefore,
\begin{equation*}
        \underset{(t,x) \to (t_0, x_0)}{\mathrm{ess}\liminf} \,  \max\left\{ \, u_t - \varphi_t + F\left(t, x, u-w, u_x - \varphi_x, u_{xx} -\varphi_{xx}\right), \, g -(u -w) \right\} \leq 0.
\end{equation*}
Combining this result with the fact that $u$ is an $L^p$-subsolution of~\eqref{eq:comparison1} then yields
\begin{equation*}
\underset{(t,x) \to (t_0, x_0)}{\mathrm{ess}\limsup} \, \max \left\{ \, \varphi_t + H\left(t,x,w, \varphi_x,  \varphi_{xx}\right), \, -w \right\} \geq 0.
\end{equation*}
This shows that $w$ is a viscosity subsolution of~\eqref{eq:proofcomparison1}.
    
By the structure condition \eqref{eq:SC}, for any test function $\varphi \in W^{1,2,p}\left(\mathcal{E}\right)$ and almost every $(t,x) \in \mathcal{E}$,
\begin{equation*}
        H\left(t,x,w, \varphi_x,  \varphi_{xx}\right) \leq \mathcal{P}_{\lambda, \Lambda}^+\left(\varphi_{xx}\right) + \omega_{1}\left( \left| w \right|\right) +\omega_{2} \left(\right|\varphi_x \left| \right).
\end{equation*}
Therefore, $w$ is also a viscosity subsolution of the problem
\begin{equation*}
\begin{cases}
    	\max \left\{ -w, \, w_t + \mathcal{P}_{\lambda, \Lambda}^+\left(w_{xx}\right) + \omega_{1}\left( \left| w \right|\right) +\omega_{2} \left(\right|w_x \left| \right) \right\} = 0 \text{ on } \mathcal{E}, \\
    	w = 0 \text{ on } \partial \mathcal{E}.
    	\end{cases}
\end{equation*}
By construction of the domain $\mathcal{E}$ and our earlier observation that $w =0$ on $\partial \mathcal{E}$, this implies that  $w$ is a viscosity subsolution of 
    \begin{equation*}
        \begin{cases}
            w_t + \mathcal{P}_{\lambda, \Lambda}^+\left(w_{xx}\right) + \omega_{1}\left( \left| w \right|\right) + \omega_{2} \left(\right|w_x \left| \right) = 0 \text{ on } \mathcal{E},\\
            w = 0 \text{ on } \partial \mathcal{E}.
        \end{cases}    
\end{equation*}
Note that $0$ is a classical solution of the above equation, and thus also a (continuous) viscosity supersolution. Theorem 4.7 in \cite{giga1991comparison} then implies that $w \leq 0$ in $\mathcal{E}$, which yields the desired  contradiction. This shows that $\mathcal{E} = \emptyset$ and, hence, that $v \geq u$ in $\bar{\mathcal{Y}}$. 
    
The proof for the case in which  $u$ is an $L^p$-supersolution, $v$ is a viscosity subsolution, and $u \geq v$ on $\partial \mathcal{Y}$ follows the same steps. \hfill $\square$
\endproof

\vspace{1em}
\begin{remark}\label{remark:comparisondegenerateellipticity}
Proposition \ref{prop:comparison} also holds when $F$ is degenerate elliptic (i.e, $\lambda = 0$ in Assumption \ref{assumption:structurecondition}) since Theorem 4.7 in \cite{giga1991comparison} applies to degenerate elliptic operators. In this case, however, the existence of $L^p$-sub- or $L^p$-supersolutions is not guaranteed.
\end{remark}

Proposition \ref{prop:comparison} and Lemma \ref{lemma:relationdifferentsolutionconcepts} yield the following corollary, which establishes the uniqueness  of $L^p$-solutions.
\begin{corollary}\label{corollary:uniqueness} 
Suppose that Assumptions \ref{assumption:structurecondition} and \ref{assumption:strictmonotonicity} hold and let $u$ and $v$ be, respectively, an $L^p$-subsolution\footnote{The same observation regarding the absence of a constraint on the functions at the frontier $\partial \mathcal{Y}$ applies as in endnote 11.} and an $L^p$-supersolution of
\begin{equation*}
        \max\left\{ \,  u_t+ F\left(t, x, u, u_x, u_{xx} \right), \,  \underset{a\in A}{\sup }\, g^a - u \right\} = 0.
\end{equation*}
If $v \geq u$ on $\partial \mathcal{Y}$, then $v \geq u$ on $\bar{\mathcal{Y}}$. 
    
If $u$ and $v$ are two $L^p$-solutions of \eqref{eq:obstacle}, then $u=v$ on $\bar{\mathcal{Y}}$.
\end{corollary}

\begin{remark} Proposition \ref{prop:comparison} and Corollary \ref{corollary:uniqueness} are also valid for non-obstacle problems, as is easily seen from the proof.\end{remark}

\subsubsection{Stability Theorem}\label{subsec:stability}

Stability theorems for both $L^p$-solutions and viscosity solutions of the Dirichlet problem (without obstacle) are well-known under a slight strengthening of the structure condition\footnote{E.g., Theorem 6.1 in \cite{crandall2000lp} gives a stability theorem when $F$ is Lipschitz continuous in $u_x$ and nonincreasing in $u$.} \eqref{eq:SC}. We establish a more general stability theorem for our setting, which is then used to prove Theorem~\ref{theorem:existenceobstacleproblemirregulardomain}. This stability theorem is valid for a large class of nonlinear parabolic obstacle problems as long as $F$ satisfies the structure condition \eqref{eq:SC}.

\begin{theorem}\label{theorem:continuityequation}
Let $\left(\mathcal{Y}^n\right)_{n \in \mathbb{N}}$ be an increasing sequence of bounded, open domains such that\footnote{$\mathcal{Y}$ does not need to be bounded.} $\bigcup_{n \in \mathbb{N}} \mathcal{Y}^n = \mathcal{Y}$. Let $\left({u}^n\right)_{n \in \mathbb{N}}$ be a sequence of functions in $W^{1,2,p}_{loc} \left( \mathcal{Y} \right) \cap \mathcal{C}^0 \left(\bar{\mathcal{Y}}\right)$ such that
\begin{enumerate}
		\item $\left(u^n\right)_{n \in \mathbb{N}}$ converges uniformly on compact subsets of $\mathcal{Y}$ and weakly in $W^{1,2,p}_{loc} \left( \mathcal{Y}\right)$ to some $u \in W^{1,2,p}_{loc} \left( \mathcal{Y}\right) \cap \mathcal{C}^0\left( \bar{\mathcal{Y}} \right)$; and
		\item for all $n \in \mathbb{N}$, $u^n = b^n$ on $\bar{\mathcal{Y}} \setminus \mathcal{Y}^n$, where the sequence $b^n$ is such that $b^n \to b$ pointwise for some function $b$ that is continuous on $\partial \mathcal{Y}$.
\end{enumerate}
Let $F$ and $F^n$, $n \in \mathbb{N}$, be operators defined on $\mathcal{Y} \times \mathbb{R} \times \mathbb{R}^d \times \mathbb{S}^d$ that satisfy Assumption \ref{assumption:structurecondition} uniformly in $n$ and Assumption \ref{assumption:convexity}. Let $g^n$, $n \in \mathbb{N}$, be continuous functions on $\bar{\mathcal{Y}}$ such that $g^n \leq b^n$ on $\bar{\mathcal{Y}} \setminus \mathcal{Y}^n$.

Suppose that for any fixed $u_{xx} \in L^p\left(\mathcal{Y}\right)$, $F^n\left(t,x, u^n, u^n_x, u_{xx} \right) \to F\left(t,x, u, u_x, u_{xx} \right)$ pointwise a.e. on $\mathcal{Y}$, and that either (i) $g^n \to g$ uniformly on compact subsets of $\mathcal{Y}$, or (ii) $g^n \uparrow g$ pointwise on $\mathcal{Y}$. If $u^n$, $n\in \mathbb{N}$, is an $L^p$-subsolution (supersolution) of 
	\begin{equation}\label{eq:limitequationn}
		\begin{cases}
			\max\left\{ u^n_t + F^n\left(t,x, u^n, u^n_x, u^n_{xx}\right) , \, g^n - u^n \right\} =0  \text{ on } \mathcal{Y}^n, \\
			u^n = b^n \text{ on } \partial \mathcal{Y}^n, \tag{Eq$^n$}
		\end{cases}
	\end{equation}
	then $u$ is an $L^p$-subsolution (supersolution) of
	\begin{equation}\label{eq:limitequation}
		\begin{cases}
			\max\left\{ u_t + F\left(t,x, u, u_x, u_{xx}\right) , \, g - u \right\} =0  \text{ on } \mathcal{Y}, \\
			u = b \text{ on } \partial \mathcal{Y}. \tag{Eq}
		\end{cases}
	\end{equation}
\end{theorem}

\vspace{0.3em}
\proof{Proof.} We first show that if, for each $n$, $u^n$ is an $L^p$-supersolution of \eqref{eq:limitequationn}, then $u$ is an $L^p$-supersolution of \eqref{eq:limitequation}.

By assumption, $u$ satisfies the regularity condition required to be an $L^p$-supersolution (requirement (i) in the definition). Moreover, since $u^n \geq b^n$ on $\partial \mathcal{Y}$ for all $n \in \mathbb{N}$,  taking limits shows that $u \geq b$ on $\partial\mathcal{Y}$ and, hence, that $u$ satisfies the second requirement for being a supersolution. 
	
We now prove the third condition, namely that 
\begin{equation*}
		 \max\left\{ u_t + F\left(t,x, u, u_x, u_{xx}\right) , \, g - u \right\}  \leq 0 \text{ a.e. on } \mathcal{Y}.
\end{equation*} 
Since $u^n$ is a $L^p$-solution of \eqref{eq:limitequationn}, (i) it is continuous, and (ii) it solves 
\begin{align*}
    & \max \left\{u^n_t + F^n(t,x,u^n, u^n_x, u^n_{xx}), \ g^n - u^n\right\} = 0 \text{ a.e. on } \mathcal{Y}^n.
\end{align*}
As a result, $g^n - u^n \leq 0 \text{ a.e. on } \mathcal{Y}^n$, and, since $u^n$ and $g^n$ are both continuous, $u^n \geq g^n$ on $\mathcal{Y}^n$. Finally, on $\bar{\mathcal{Y}} \setminus \mathcal{Y}_n$, $u^n = b^n \geq g^n$ by Assumption \ref{assumption:gregularity}.(ii). Thus, $u^n \geq g^n$ on $\bar{\mathcal{Y}}$. So, for all $n \in \mathbb{N}$, $u^n \geq g^n$ on $\bar{\mathcal{Y}}$. Passing to the limit, $u \geq g$ on $\bar{\mathcal{Y}}$. 
	
There remains to show that
\begin{equation*}
		u_t + F\left(t,x, u, u_x, u_{xx}\right) \leq  0 \text{ a.e. on } {\mathcal{Y}},
\end{equation*}
which will hold if \begin{equation*}
		u_t + F\left(t,x, u, u_x, u_{xx}\right) \leq  0 \text{ a.e. on } {\mathcal{O}}
\end{equation*}
for any strict bounded subset $\mathcal{O}$ of\footnote{We say that $\mathcal{O}$ is a strict subset of $\mathcal{Y}$ if the closure of $\mathcal{O}$ lies in the interior of $\mathcal{Y}$.} $\mathcal{Y}$.

Thus, consider any  strict bounded subset $\mathcal{O}$ of $\mathcal{Y}$. Since $\mathcal{Y} = \bigcup_{n \in \mathbb{N}} \mathcal{Y}^n$, compactness implies that there exists some integer $N\geq 1$ such that $\mathcal{Y}^N$ contains $\mathcal{O}$. Moreover, since the sequence $\{\mathcal{Y}^n\}_{n\geq 1}$ is increasing, $\mathcal{O} \subset \mathcal{Y}^n$ for all $n \geq N$. Without loss of generality, we assume that $N = 1$. For all $m \in \mathbb{N}$, define
\begin{equation*}
		\underline{H}^m\left( t,x, v_t, v_{xx} \right) = \underset{k\geq m}{\inf } \, \left\{v_t + F^k \left(t,x, u^k, u^k_x, v_{xx}\right)\right\}.
\end{equation*}
Since $u^n_t + F^n \left(t,x, u^n, u^n_x, u^n_{xx}\right) \leq 0$, we have for $n\geq m$ and almost every $(t,x) \in \mathcal{O}$,
\begin{equation*}
		\underline{H}^m\left( t,x, u^n_t, u^n_{xx} \right) \leq 0.
\end{equation*}
Since $F^n$ satisfies Assumptions \ref{assumption:structurecondition} and \ref{assumption:convexity} for each $n$, we deduce that for all $m\in \mathbb{N}$ and all $v_t \in \mathbb{R}$, $\underline{H}^m(t,x, v_t, \cdot )$ is Lipschitz continuous in $v_{xx}$ a.e. on $\mathcal{O}$ and
\begin{equation*}
\lambda I_{d} \leq_{\mathbb{S}^d} D_{v_{xx}} \underline{H}^m(t,x, v_t, \cdot) \leq_{\mathbb{S}^d} {\Lambda} I_d.
\end{equation*}
Moreover, for any $v_{xx} \in \mathbb{S}^d$, the structure condition \eqref{eq:SC} implies that $\underline{H}^m$ is bounded in $L^p\left( \mathcal{O}\right)$ because the functions $u^n$, $n \in \mathbb{N}$, are uniformly bounded in $W^{1,2,p}_{loc}\left( \mathcal{Y}\right)$ as part of a weakly convergent sequence. Finally, by assumption, $u^n \rightharpoonup u$ in $W^{1,2,p}\left( \mathcal{O} \right)$. So, $u$ and its weak derivatives are bounded in $L^p$ by the uniform boundedness principle, and hence, $\sup_n \underline{H}^m(t,x, u_t^n, u_{xx}^n)$ is bounded in $L^p$ by the  structure condition \eqref{eq:SC} and Assumption \ref{assumption:growthconditions}. Theorem 4.2.6 (b) in \cite{krylov2018sobolev} then implies that, for all\footnote{In order to apply Theorem 4.2.6 in \cite{krylov2018sobolev} to the operator $\underline{H}^m$, we use the fact that it is a $\lambda$-nondegenerate $\mathcal{L}$-type operator, which follows from parts (c) and (d) of Lemma 4.2.4 in \cite{krylov2018sobolev}: part (c) shows that $F^n$ has this property, and part (d) shows that the supremum over $n$ also inherits the property.} $m \in \mathbb{N}$,
\begin{equation*}
    \underline{H}^m\left( t,x, u_t, u_{xx} \right) = \underset{k\geq m}{\inf } \left\{u_t + F^k \left(t,x, u^k, u^k_x, u_{xx}\right) \right\} \leq \underset{n\to \infty}{\limsup}\, \underline{H}^m\left( t,x, u^n_t, u^n_{xx} \right) \leq 0 \text{ a.e. on } \mathcal{O}.
\end{equation*}
Moreover, $ u_t+ F^n \left(t,x, u^n, u^n_x, u_{xx}\right) \to u_t + F \left(t,x, u, u_x, u_{xx}\right)$ pointwise a.e. on $\mathcal{Y}$ by assumption. Letting $m\to \infty$, we obtain 
\begin{equation*}u_t + F\left(t,x, u, u_x, u_{xx}\right) \leq 0 \text{ a.e. on } \mathcal{O}.
\end{equation*}
This shows that $u$ is an $L^p$-supersolution of \eqref{eq:limitequation}.
	
Next, we prove that if, for each $n$, ${u^n}$ is an $L^p$-subsolution of \eqref{eq:limitequationn}, then $u$ is an $L^p$-subsolution of \eqref{eq:limitequation}.

We first note that in this case, we have  $u^n \leq b^n$ on $\partial \mathcal{Y}$ for all $n \in \mathbb{N}$, . Taking limits, $u \leq b$ on $\partial\mathcal{Y}$.

Next, we show that
\begin{equation*}\max\left\{u_t + F\left(t,x, u, u_x, u_{xx}\right), \, g-u\right\} \geq  0 \text{ a.e. on } \mathcal{Y},\end{equation*}
which holds if and only if 
\begin{equation*}\max\left\{u_t + F\left(t,x, u, u_x, u_{xx}\right), \, g-u\right\} \geq  0 \text{ a.e. on } \mathcal{O},\end{equation*}
for any strict bounded subset $\mathcal{O}$ of $\mathcal{Y}$. 

Let $\mathcal{O}$ be a strict bounded subset of $\mathcal{Y}$. Since $\mathcal{Y} = \bigcup_{n \in \mathbb{N}} \mathcal{Y}^n$ and the sequence of subsets $\mathcal{Y}^n$ is increasing, there exists $N \in \mathbb{N}$ such that $\mathcal{O} \subset \mathcal{Y}^n$, for all $n \geq N$. Without loss of generality, assume that $N = 1$. It is then enough to show that, on $\mathcal{E} = \left\{ (t,x) \in \mathcal{O} \, : \, u(t,x) > g(t,x) \right\}$, we have
\begin{equation*}u_t + F\left(t,x, u, u_x, u_{xx} \right) \geq 0 \text{ a.e. on } \mathcal{O}.\end{equation*}
If $\mathcal{E} = \emptyset$, there is nothing left to prove. Suppose instead that $\mathcal{E}$ is nonempty, and let $k$ be large enough so that
\begin{equation*}\mathcal{E}^k = \left\{ (t,x) \in \mathcal{O} \, : \, u(t,x) > g(t,x) +\frac{1}{k} \right\} \neq \emptyset.\end{equation*}
Then, there exists $N^k \in \mathbb{N}$ such that $g^n < u^n$ on $\mathcal{E}^k$ for all $n \geq N^k$. To see this, choose $N^k$ such that $\left| u^n - u \right| < \frac{1}{3k}$ on $\mathcal{E}^k$ (which exists since $u^n \to u$ uniformly). Then recall that either (i) $g^n \uparrow g$ pointwise, or (ii) $g^n \to g$ uniformly. So, either (i) $g^n \leq g < u - \frac{1}{k} < u - \frac{1}{3k} < u^n$ on $\mathcal{E}^k$ and we are done, or (ii) $g^n$ converge uniformly to $g$. In that case, choose $\tilde{N}^k$ such that $\left| g^n - g \right| < \frac{1}{3k}$ for all $(t,x) \in \mathcal{O}$, and again, it follows that $g^n < g + \frac{1}{3k} < u - \frac{1}{3k} < u^n$ on $\mathcal{E}^k$ for all $n \geq \max\left\{ {N}^k, \tilde{N}^k  \right\}$. 
 
Since $g^n < u^n$ on $\mathcal{E}^k$ for $n > N^k$, we conclude from \eqref{eq:limitequationn} that, for all $n \geq N^k$, 
\begin{equation*}u^n_t + F^n\left(t,x, u^n, u^n_x, u^n_{xx} \right) = 0 \text{ a.e. on } \mathcal{E}^k.	\end{equation*}
Similarly to the first part of the proof, define for all $m \geq N^k$
\begin{equation*}\bar{H}^m\left( t,x, v_t, v_{xx} \right) = \underset{l\geq m}{\sup } \, \left\{ v_t + F^l \left(t,x, u^l, u^l_x, v_{xx}\right)\right\}.\end{equation*}
We note once more that, for $n\geq m$ and for almost every $(t,x)\in \mathcal{E}^k$,
\begin{equation*}\bar{H}^m\left( t,x, u^n_t, u^n_{xx} \right) \geq 0.\end{equation*}
Moreover, as above, for all $m\in \mathbb{N}$ and all $v_t \in \mathbb{R}$, $\bar{H}^m(t,x, v_t, \cdot )$ is Lipschitz continuous in $v_{xx}$ a.e. on $\mathcal{E}^k$ with 
\begin{equation*}
\lambda I_d \leq_{\mathbb{S}^d} D_{v_{xx}} \bar{H}^m(t,x, v_t, \cdot) \leq_{\mathbb{S}^d} {\Lambda} I_d.
\end{equation*}
Also, for any $v_{xx} \in \mathbb{S}^d$, by the structure condition \eqref{eq:SC}, $\bar{H}^m$ is bounded in $L^p\left( \mathcal{O}\right)$ because the functions $u^n$, $n \in \mathbb{N}$, are  uniformly bounded in $W^{1,2,p}_{loc}\left( \mathcal{Y}\right)$ as they form a weakly convergent sequence. Finally, by assumption, $u^n \rightharpoonup u$ in $W^{1,2,p}\left( \mathcal{E}^k \right)$. So, $u$ and its weak derivatives are bounded in $L^p$ by the uniform boundedness principle, and hence, $\inf_n \bar{H}^m(t,x, u_t^n, u_{xx}^n)$ is bounded in $L^p$ by the  structure condition \eqref{eq:SC} and Assumption \ref{assumption:growthconditions}. Theorem 4.2.6 (c) in \cite{krylov2018sobolev} then implies that, for all\footnote{Theorem 4.2.6 in \cite{krylov2018sobolev} also applies to the operator $\bar{H}^m$ by the observation made in the previous endnote.} $m \in \mathbb{N}$,
\begin{equation*}
\bar{H}^m\left( t,x, u_t, u_{xx} \right) = \underset{l\geq m}{\sup } \left\{ u_t + F^l \left(t,x, u^l, u^l_x, u_{xx}\right) \right\} \geq \underset{n\to \infty}{\liminf} \, \bar{H}^m\left( t,x, u^n_t, u^n_{xx} \right)\geq 0 \text{ a.e. on } \mathcal{E}^k.
\end{equation*}
By assumption, $u_t + F^n \left(t,x, u^n,u^n_x, u_{xx}\right) \to u_t + F \left(t,x, u, u_x, u_{xx}\right)$ pointwise a.e. on $\mathcal{Y}$. Letting $m\to \infty$, we obtain 
\begin{equation*}u_t + F\left(t,x, u, u_x, u_{xx}\right) \geq 0 \text{ a.e. on } \mathcal{E}^k.\end{equation*}
	Since $\bigcup_{k\in \mathbb{N}} \mathcal{E}^k = \mathcal{E}$, it follows that
	\begin{equation*}
		u_t + F\left(t,x, u, u_x, u_{xx} \right)  \geq 0 \text{ a.e. on } \mathcal{E}.
	\end{equation*}
	This shows that $u$ is an $L^p$-subsolution of \eqref{eq:limitequation}. \hfill$\square$

\begin{remark}
    Theorem \ref{theorem:continuityequation} is also valid for non-obstacle problems, as can be seen from the proof. Moreover, it also holds for viscosity solutions. To see this, one can invoke Theorem 6.1 in \cite{crandall2000lp} instead of Theorem 4.2.6 in \cite{krylov2018sobolev} in the proof.
\end{remark}

\subsection{Existence, uniqueness, and $W^{1,2,p}$-estimate when $\left|A\right| < \infty$}\label{sec:finiteA}

Before considering the general case, we derive a more restrictive version of Theorem~1 when the obstacle is the maximum of finitely many $W^{1,2,p}$-obstacles and the following strengthening Assumptions 1 and 2 is imposed.

\begin{assumptionp}{1'}\label{assumption:smoothdomain}
    $\mathcal{X}$ is an $\mathcal{C}^{1, Lip}$ open bounded subset of $\mathbb{R}^d$.
\end{assumptionp}

\begin{assumptionp}{2'}\label{assumption:bregularity}
	${b}: \bar{\mathcal{Y}} \to \mathbb{R}$ is in $W^{1,2,p}\left( \bar{\mathcal{Y}} \right)$.
\end{assumptionp}

\begin{lemma}\label{theorem:existenceobstacleproblemfiniteA}
Suppose that Assumptions \ref{assumption:smoothdomain}, \ref{assumption:bregularity}, and  \ref{assumption:gregularity}--\ref{assumption:growthconditions} hold for some\footnote{The parameter $p$ appears in Assumptions~3 and 7.} $p \in (d+2, \infty)$. If $\left|A\right| < \infty$, then \eqref{eq:obstacle} has an $L^p$-solution $u$. Moreover, $u \in W^{1,2,p}\left(\mathcal{Y}\right)$ and there exists $C = C \left( d, p, \lambda, \Lambda, R, \eta_F, diam(\mathcal{X}), T, L^{1,Lip}(\mathcal{X}) \right) \in \mathbb{R}_+$ such that
\begin{equation}\label{eq:W12pestimatesmoothdomainfiniteA}
\left\| u \right\|_{W^{1,2,p}\left(\mathcal{Y}\right)} \leq C \left( 1 + \underset{a \in A}{\max}\left\| g^a \right\|_{W^{1,2,p} \left( \mathcal{Y}\right)} + \left\| b \right\|_{W^{1,2,p}\left(\mathcal{Y}\right)} + \left\| G \right\|_{L^p\left(\mathcal{Y}\right)} \right).
\end{equation}

If Assumption 8 also holds, the solution $u$ is unique. 
\end{lemma}

The proof builds on the approximation argument proposed by \cite{byun2018nondivergence} and \cite{byun2022w}, extending it to fully nonlinear parabolic obstacle problems. \cite{byun2018nondivergence} approximate the obstacle problem by a Dirichlet problem similar to \eqref{eq:Dirichlet}. They then use the Schauder fixed-point theorem to show the existence of $L^p$-solutions whose norms are bounded by a constant independent of $\epsilon$. Instead, we appeal to known existence results and estimates on the regularity of solutions of the Dirichlet problem, e.g., Theorem 12.1.7, 15.1.3, and 15.1.4 in \cite{krylov2018sobolev} that allow for an arbitrary continuous dependence of the operator on the value of the solution $u$. As a result, we can sidestep the fixed-point argument in \cite{byun2018nondivergence}, streamlining the proof of the existence of a solution and $W^{1,2,p}$-estimate when the obstacle is in $W^{1,2,p}$. We then show by induction that the existence and regularity results for a single $W^{1,2,p}$-obstacle carry to the pointwise maximum of a finite number of $W^{1,2,p}$-obstacles.

\vspace{0.3em}
\proof{Proof.}
Without loss of generality, assume that $A = \left\{1,\dots, I\right\}$, for some integer $I\geq 1$. The proof proceeds by induction on $I$. For all $I \geq 1$, define property $\bf{P(I)}$ as follows:
\vspace{0.5em}
\begin{quote}
    If $g = \underset{a \in \left\{ 1,\dots, I \right\}}{\max }\, g^a$ with $g^a \in W^{1,2,p}\left(\mathcal{Y}\right)$ for all $a \in \left\{ 1,\dots, I \right\}$ and the operator $F$, domain $\mathcal{Y}$, and boundary data $b$ in \eqref{eq:obstacle} satisfy Assumptions \ref{assumption:smoothdomain}, \ref{assumption:bregularity}, and  \ref{assumption:gregularity}--\ref{assumption:growthconditions} for some $p \in (d+2, \infty)$, then \eqref{eq:obstacle} has a $L^p$-solution $u \in W^{1,2,p}\left(\mathcal{Y}\right)$. Moreover, $u$ satisfies the estimate \eqref{eq:W12pestimatesmoothdomainfiniteA}.
\end{quote}
\vspace{0.5em}

We start by proving the  base case $\bf{P(1)}$. Thus, let $g^1 \in W^{1,2,p}\left(\mathcal{Y}\right)$ be such that  $g^1 \leq b$ on $\partial \mathcal{Y}$. Let  $h^1(t,x) = - g^1_t - F\left(t, x, g^1, g^1_x, g^1_{xx} \right)$. For each $\epsilon >0$, let $\Phi_{\epsilon} \in \mathcal{C}^{\infty}\left(\mathbb{R}\right)$ be a nondecreasing function such that $\Phi_{\epsilon}(a) = 0$ if $a \leq 0$, and $\Phi_{\epsilon}(a) = 1$ if $a \geq \epsilon$. In particular, for all $a \in \mathbb{R}$, $\Phi_{\epsilon}(a) \in [0,1]$. Finally, consider a sequence  $\left(\epsilon_n\right)_{n \in \mathcal{N}} \subseteq \mathbb{R}_{++}$ such that $\epsilon_n \to 0$ as $n \to \infty$, and consider the following auxiliary nonlinear Dirichlet problem for each $n \in \mathbb{N}$:
\begin{equation}\label{eq:Dirichlet}
    \begin{cases}
        u^{n}_t + {F}\left(t, x, u^{n}, u^{n}_x, u^{n}_{xx} \right) = h^1(t,x)^+ \Phi_{\epsilon_n} \left( u^n - g^1 \right) - h^1(t,x)^+  \text{ on } \mathcal{Y}, \\
        u^{n} = b \text{ on } \partial\mathcal{Y}.
    \end{cases} \tag{D$^n$}
\end{equation}
Here $h^1(t,x)^+ = \max\left\{ h^1(t,x), 0\right\}$ stands for the nonnegative part of $h^1$. By Assumptions \ref{assumption:structurecondition} and \ref{assumption:growthconditions}, $h^1(t,x)$ is in ${L}^{p}\left(\mathcal{Y}\right)$ and
\begin{equation}\label{eq:boundonh}
	\left\| h^1 \right\|_{L^{p}(\mathcal{Y})} \leq  C^{h^1} \left( 1+ \left\| g^1 \right\|_{W^{1,2,p}(\mathcal{Y})} + \left\| G \right\|_{L^{p}(\mathcal{Y})} \right)
\end{equation}
where $C^{h^1} = C\left(d,p, \lambda, \Lambda, {R}\right)\in \mathbb{R}_+$. By Theorem 15.1.4 in \cite{krylov2018sobolev}, for each $n\in \mathbb{N}$, there exists a solution $u^n \in W^{1,2,p}\left( \mathcal{Y} \right)$ of \eqref{eq:Dirichlet}. Moreover, by Theorem 15.1.3 in the same work,
\begin{equation*}
    \left\| u^n \right\|_{W^{1,2,p}\left(\mathcal{Y}\right)} \leq C \left( 1 +\left\| h^1 \right\|_{L^{p}(\mathcal{Y})} +  \left\| G \right\|_{L^p\left(\mathcal{Y}\right)} + \left\| b \right\|_{W^{1,2,p}\left(\mathcal{Y}\right)} + \left\| u^n \right\|_{L^\infty\left(\mathcal{Y}\right)} \right)
\end{equation*}
where $C = C \left( d, p, \lambda, \Lambda, R, diam(\mathcal{X}), T, L^{1,Lip}(\mathcal{X}) \right) \in \mathbb{R}_+$. By Lemma 12.1.9 in \cite{krylov2018sobolev}, 
\begin{equation*}
     \left\| u^n \right\|_{L^{\infty}\left(\mathcal{Y}\right)} \leq C  \left( 1 + \left\| h^1 \right\|_{L^{p}(\mathcal{Y})} + \left\| G \right\|_{L^p\left(\mathcal{Y}\right)}+ \left\| b \right\|_{L^{\infty}\left(\partial \mathcal{Y}\right)} \right),
\end{equation*}
where $C = C \left( d, p, \lambda, \Lambda, R, \eta_F, diam(\mathcal{X}), T \right) \in \mathbb{R}_+$. Finally, $\left\| b \right\|_{L^\infty\left( \partial \mathcal{Y}\right)} \leq $$\left\| b \right\|_{L^\infty\left(\mathcal{Y}\right)} \leq C \left\| b\right\|_{W^{1,2,p}\left( \mathcal{Y} \right)}$, where $C = C(d, p, diam(\mathcal{X}), T, L^{0,Lip}(\mathcal{X}))$, by the Morrey-Sobolev embedding theorem. Combining the inequalities above, we obtain
\begin{equation*}
	\left\| u^n \right\|_{W^{1,2,p}\left(\mathcal{Y}\right)} \leq C \left( 1 + \left\| g^1 \right\|_{W^{1,2,p} \left( \mathcal{Y}\right)} + \left\| b \right\|_{W^{1,2,p}\left(\mathcal{Y}\right)} + \left\| G \right\|_{L^p\left(\mathcal{Y}\right)} \right)
\end{equation*}
where $C = C \left( d, p, \lambda, \Lambda, R, \eta_F, diam(\mathcal{X}), T, L^{0,Lip}(\mathcal{X}) \right)\in \mathbb{R}_+$ is independent of $n$.
    
$W^{1,2,p}(\mathcal{Y})$ is separable and reflexive.\footnote{The space ${L}^p$ is uniformly convex for $p \in (1, \infty)$, hence so is $W^{1,2,p}\left(\mathcal{Y}\right)$. Therefore, this latter space is reflexive by Theorem 1.40 in \cite{demengel2012functional}.} Therefore, its closed bounded subsets are weakly sequentially compact by Theorem 1.32 in \cite{demengel2012functional}. Moreover $W^{1,2,p}\left(\mathcal{Y}\right)$ is compactly embedded in $\mathcal{C}^{0}\left(\bar{\mathcal{Y}}\right)$ by the Rellich-Kondrachov theorem (see Theorem 2.84 in \cite{demengel2012functional}). Therefore, there exists a function $u^{(1)} \in W^{1,2,p}\left(\mathcal{Y}\right) \cap \mathcal{C}^{0}\left(\bar{\mathcal{Y}}\right)$ and a subsequence $\left( u^{n_j}\right)_{j \in \mathbb{N}} \subseteq \left( u^{n}\right)_{n\in \mathbb{N}}$ such that
\begin{equation*}
    \begin{cases}
        u^{n_j}  \rightharpoonup u^{(1)} \text{ in } W^{1,2,p}\left( \mathcal{Y} \right),\\
        u^{n_j} \to u^{(1)} \text{ in } \mathcal{C}^{0}\left( \bar{\mathcal{Y}} \right),
    \end{cases}
\end{equation*}
as $j \to \infty$. Moreover, $u^{(1)}$ satisfies the estimate \eqref{eq:W12pestimatesmoothdomainfiniteA}:
\begin{equation*}
    \begin{aligned}
    \left\| u^{(1)} \right\|_{W^{1,2,p}\left(\mathcal{Y}\right)} & \leq \underset{j \to \infty}{\liminf}\, \left\| u^{n_j} \right\|_{W^{1,2,p}\left(\mathcal{Y}\right)} \\
    & \leq C \left( 1 + \left\| g^1 \right\|_{W^{1,2,p} \left( \mathcal{Y}\right)} + \left\| b \right\|_{W^{1,2,p}\left(\mathcal{Y}\right)} + \left\| G \right\|_{L^p\left(\mathcal{Y}\right)} \right).
    \end{aligned}
\end{equation*}
There remains to show that $u^{(1)} \in W^{1,2,p}\left(\mathcal{Y}\right)$ is an $L^p$-solution of \eqref{eq:obstacle}. 

First, for all $j \in \mathbb{N}$, $u^{n_j} = b$ on $\partial \mathcal{Y}$, and $\left(u^{n_j}\right)_{j \in \mathbb{N}}$ converges uniformly to $u$ on $\bar{\mathcal{Y}}$. So, $u^{(1)} = b$ on $\partial{\mathcal{Y}}$.
	
Second, we prove that
\begin{equation*}
    u^{(1)}_t + F\left(t,x, u^{(1)}, u^{(1)}_x, u^{(1)}_{xx} \right) \leq 0 \text{ a.e. on } \mathcal{Y}.
\end{equation*}
For all $j \in \mathbb{N}$ and a.e. $(t,x) \in \mathcal{Y}$,
\begin{equation*}
	u^{n_j}_t + {F}\left(t, x, u^{{n_j}}, u_x^{{n_j}}, u_{xx}^{{n_j}}\right) = {h^1}(t,x)^+ \Phi_{\epsilon_{n_j}} \left( u^{{n_j}} -g^1 \right) - {h^1}(t,x)^+ \leq 0,
\end{equation*}
So, for all $j \in \mathbb{N}$, $u^{n_j}$ is an $L^p$-supersolution of
\begin{equation*}
	u^{n_j}_t + {F}\left(t, x, u^{{n_j}}, u_x^{{n_j}}, u_{xx}^{{n_j}}\right) = 0.
\end{equation*}
Theorem \ref{theorem:continuityequation} then implies that $u^{(1)}$ is an $L^p$-supersolution of
\begin{equation*}
	u^{(1)}_t+ F\left(t,x, u^{(1)}, u^{(1)}_x, u^{(1)}_{xx} \right) = 0 \text{ on } \mathcal{Y}.
\end{equation*}
To see this, To see this, recall that the subsequence $\left( u^{n_j}\right)_{j \in \mathbb{N}}$ converges uniformly on compact subsets of $\bar{\mathcal{Y}}$ and weakly in $W^{1,2,p}_{loc}(\mathcal{Y})$ to $u^{(1)} \in W^{1,2,p}_{loc}(\mathcal{Y}) \cap \mathcal{C}^0(\bar{\mathcal{Y}})$. As a result, passing to a further subsequence if necessary, we can assume that $u^{n_j}_{x}$ also converges to $u^{(1)}_x$ uniformly on compact subsets $\mathcal{Y}$ by a Morrey-Sobolev embedding type theorem (Lemma 3.3, page 80 in \cite{ladyzenskaja1968linear}). Hence, by the structure condition \eqref{eq:SC}, for all fixed $u_{xx}$, ${F}\left(t, x, u^{n_j}, u^{n_j}_x, u_{xx} \right) \to F\left(t,x, u^{(1)}, u^{(1)}_x, u_{xx} \right)$ pointwise a.e.. Therefore, Theorem \ref{theorem:continuityequation} applies and yields the desired result.

Next, we show that $u^{(1)} \geq g^1$ a.e. on $\mathcal{Y}$. For each $j \in \mathbb{N}$, define
\begin{equation*}
	\mathcal{V}_{n_j} = \left\{ (t,x) \in \mathcal{Y} \, : \, g^1(t,x) > u^{{n_j}}(t,x) \right\}.
\end{equation*}
We will show that  $\mathcal{V}_{n_j}$ is empty for each $j \in \mathbb{N}$. Suppose by way of contradiction that $\mathcal{V}_{n_j} \neq \emptyset$ for some $j \in \mathbb{N}$. Since both $u^{n_j}$ and $g^1$ are continuous on $\mathcal{Y}$, $\mathcal{V}_{n_j}$ is open. Moreover, since $b \geq g^1$, $u^{n_j} = g^1$ on $\partial \mathcal{V}_{n_j}$. On $\mathcal{V}_{n_j}$, $\Phi_{\epsilon_{n_j}}\left(u^{{n_j}}(t,x) - g(t,x)\right) = 0$. So $u^{{n_j}} \in W^{1,2,p}\left( \mathcal{Y} \right)$ is an $L^p$-solution of
\begin{equation*}
	\begin{cases}
    	u^{n_j}_t + F\left(t, x, u^{{n_j}}, u_x^{{n_j}}, u_{xx}^{{n_j}}\right) = - {h^1}(t,x)^+ \text{ on } \mathcal{V}_{n_j}, \\
    	u^{n_j} = g^1 \text{ on } \partial \mathcal{V}_{n_j}.
     \end{cases}
\end{equation*}
Since $ - {h^1}^+ \leq -{h^1}$, $u^{{n_j}}$ is an $L^p$-solution of
\begin{equation*}
	\begin{cases}
		u^{n_j}_t + {F}\left(t, x, u^{n_j}, u_x^{n_j}, u_{xx}^{n_j}\right) \leq - {h^1}(t,x) \text{ on } \mathcal{V}_{n_j} \\
		u^{n_j} = g^1 \text{ on } \partial \mathcal{V}_{n_j}.
	\end{cases}
\end{equation*}
By definition of ${h^1}$, $g^1 \in W^{1,2,p} \left( \mathcal{Y} \right)$ is an $L^p$-solution of
\begin{equation*}
	\begin{cases}
    	g^1_t + {F}\left(t, x, g^1, g^1_x, g^1_{xx} \right) = - {h^1}(t,x) \text{ on } \mathcal{V}_{n_j} \\
    	g^1 = g^1 \text{ on } \partial \mathcal{V}_{n_j}.
    \end{cases}
\end{equation*}
Therefore, by our comparison principle (Proposition \ref{prop:comparison}), $g^1 \leq u^{n_j}$ on $\bar{\mathcal{V}}_{n_j}$, which yields the desired  contradiction. We have thus showed that,for all $j\in \mathbb{N}$, $\mathcal{V}_{n_j} = \emptyset$, i.e., $u^{n_j} \geq g^1$ on $\bar{\mathcal{Y}}$. Taking limits, as $\left(u^{n_j}\right)_{j \in \mathbb{N}}$ converges uniformly to $u^{(1)}$ on $\bar{\mathcal{Y}}$, $u^{(1)} \geq g^1$ on $\mathcal{Y}$.
	
We have shown so far that $u^{(1)} = b$ on $\partial \mathcal{Y}$ and that
	\[\max\left\{ u^{(1)}_t + F\left(t,x,u^{(1)}, u^{(1)}_x, u^{(1)}_{xx} \right), \, g^1-u^{(1)}\right\} \leq 0 \text{ a.e. on } \mathcal{Y},\]
i.e., that $u$ is an $L^p$-supersolution of \eqref{eq:obstacle}, with $u = b$ on $\partial \mathcal{Y}$. To conclude, there remains to show that
\begin{equation*}
	\max\left\{ u^{(1)}_t + F\left(t,x,u^{(1)}, u^{(1)}_x, u^{(1)}_{xx} \right), \, g^1-u^{(1)}\right\} \geq 0 \text{ a.e. on } \mathcal{Y}.
\end{equation*}
We do so by showing that, on the open set $\mathcal{U} = \left\{ (t,x) \in \mathcal{Y} \, : \, u^{(1)}(t,x) > g^1(t,x) \right\}$, we have
\begin{equation*}
	u^{(1)}_t + F\left(t,x, u^{(1)}, u^{(1)}_x, u^{(1)}_{xx} \right) \geq 0 \text{ a.e.}.
\end{equation*}
This follows from Theorem \ref{theorem:continuityequation}. To see this, observe that each $u^{{n_j}}$ solves the associated Dirichlet problem \eqref{eq:Dirichlet}. But, by definition of $\Phi_{\epsilon}$, we have $\Phi_{\epsilon_{n_j}} \left( u^{{n_j}}(t,x) - g(t,x) \right) \to 1$ pointwise a.e. on $\mathcal{U}$ as $j \to \infty$, and, hence, the right-hand side of the partial differential equation \eqref{eq:Dirichlet} converges pointwise a.e. to zero. Moreover, for all $(t_0, x_0) \in \mathcal{U}$ and $D>0$ such that $C_D(t_0, x_0) \subset \mathcal{U}$, the subsequence $\left( u^{n_j}\right)_{j \in \mathbb{N}}$ converges uniformly on compact subsets of $\bar{\mathcal{Y}}$, hence on compact subset of $C_D(t_0, x_0) \subset \mathcal{U}$, and weakly in $W^{1,2,p}_{loc}(C_D(t_0, x_0))$ to $u^{(1)} \in W^{1,2,p}(\mathcal{C}_D(t_0,x_0)) \subset  W^{1,2,p}_{loc}({C}_D(t_0,x_0)) \cap \mathcal{C}^0(\bar{C}_D(t_0, x_0))$. As a result, passing to a further subsequence if necessary, we can assume that $u^{n_j}_{x}$ also converges to $u^{(1)}_x$ uniformly on compact subsets ${C}_D(t_0,x_0)$ by a Morrey-Sobolev embedding type theorem (Lemma 3.3, page 80 in \cite{ladyzenskaja1968linear}). Hence, by the structure condition \eqref{eq:SC}, for all fixed $u_{xx}$, ${F}\left(t, x, u^{n_j}, u^{n_j}_x, u_{xx} \right) - h^1(t,x)^+ \Phi_{\epsilon_n} \left( u^{n_j} - g^1 \right) + h^1(t,x)^+ \to F\left(t,x, u^{(1)}, u^{(1)}_x, u_{xx} \right)$ pointwise a.e.. Therefore, Theorem \ref{theorem:continuityequation} applies, and, for all $(t_0, x_0) \in \mathcal{U}$ and $D>0$ such that $C_D(t_0, x_0) \subset \mathcal{U}$, $u$ is an $L^p$-subsolution of
\begin{equation*}
	u^{(1)}_t + F\left(t,x, u^{(1)}, u^{(1)}_x, u^{(1)}_{xx} \right) = 0 \text{ on } C_D(t_0, x_0).
\end{equation*}
Since $(t_0, x_0) \in \mathcal{U}$ was arbitrary, we obtain
\begin{equation*}
	\max\left\{ u^{(1)}_t + F\left(t,x,u^{(1)}, u^{(1)}_x, u^{(1)}_{xx} \right), \, g^1-u^{(1)}\right\}  \geq 0 \text{ a.e. on } \mathcal{U}.
\end{equation*}
Since $g^1\geq u^{(1)}$ on the complement of $\mathcal{U}$, we conclude that
\begin{equation*}
	\max\left\{ u^{(1)}_t + F\left(t,x,u^{(1)}, u^{(1)}_x, u^{(1)}_{xx} \right), \, g^1-u^{(1)}\right\} \geq 0 \text{ a.e. on } \mathcal{Y}.
\end{equation*}
and, hence, that $u^{(1)}$ is an $L^p$-solution of \eqref{eq:obstacle} with $I=1$. This concludes the base case of our induction argument.

Next, suppose that $\bf{P(I)}$ holds for some $I \in \mathbb{N}$. We will show that it also holds for $I+1$. Let $g^{I+1} \in W^{1,2,p}\left( \mathcal{Y}\right)$ with $g^{I+1} \leq b$ on $\partial \mathcal{Y}$. Define $h^{I+1} = - g^{I+1}_t - F\left(t, x, g^{I+1}, g^{I+1}_x, g^{I+1}_{xx} \right) $; and consider the following sequence of auxiliary problems. For all $n \in \mathbb{N}$,
\begin{equation}\label{eq:Dirichletinduction}
    \begin{cases}
        \max \bigg\{ u^n_t + F\left(t,x,u^n, u^n_x, u^n_{xx}\right) -h^{I+1}(t,x)^+ \Phi_{\epsilon_{n}} \left( u^{{n}} -g^{I+1} \right) + h^{I+1}(t,x)^+, \\
        \qquad \qquad \qquad \qquad \qquad \qquad \underset{i\in \{1,\dots, I\}}{\max}\,g^i - u^n \bigg\} = 0 \text{ on } \mathcal{Y}, \\
        u^n =b \text{ on } \partial \mathcal{Y}.
    \end{cases}
\end{equation}
Observe that, for each $n\in \mathbb{N}$, the obstacle $g$, boundary data $b$, domain $\mathcal{Y}$, and operator $F\left(t,x,\cdot, \cdot, \cdot \right) -h^{I+1}(t,x)^+ \Phi_{\epsilon_{n}} \left( \cdot -g^{I+1} \right) + h^{I+1}(t,x)^+$ in \eqref{eq:Dirichletinduction} satisfy the assumptions of our induction hypothesis. So, by our induction hypothesis, for all $n \in \mathbb{N}$, \eqref{eq:Dirichletinduction} has an $L^p$-solution $u^n \in W^{1,2,p}\left(\mathcal{Y}\right)$, with 
\begin{equation*}
    \left\| u^n \right\|_{W^{1,2,p}\left( \mathcal{Y} \right)} \leq C \left( 1 + \underset{i \in \{ 1, \dots, I\}}{\max} \, \left\| g \right\|_{W^{1,2,p} \left( \mathcal{Y}\right)} + \left\| b \right\|_{W^{1,2,p}\left(\mathcal{Y}\right)} + \left\| G \right\|_{L^p\left(\mathcal{Y}\right)} + \left\| g^{I+1} \right\|_{W^{1,2,p} \left( \mathcal{Y}\right)} \right),
\end{equation*}
where $C = C \left( d, p, \lambda, \Lambda, \eta_F, diam(\mathcal{X}), T, L^{1,Lip}(\mathcal{X}) \right)\in \mathbb{R}_+$ is independent of $\epsilon_n$, since
\begin{equation*}
    \left\| h^{I+1} \right\|_{L^{p}(\mathcal{Y})} \leq  C^{h} \left( 1+ \left\| g^{I+1} \right\|_{W^{1,2,p}(\mathcal{Y})} + \left\| G \right\|_{L^{p}(\mathcal{Y})} \right).    
\end{equation*}
Proceeding as in the base case $I=1$, we see that the sequence $\left(u^n\right)_{n\in\mathbb{N}}$ has a subsequence that converges weakly in $W^{1,2,p}\left(\mathcal{Y}\right)$ and strongly in $C^0\left(\bar{\mathcal{Y}}\right)$ to some $u^{(I+1)} \in W^{1,2,p}\left(\mathcal{Y}\right)$. Next, we show that $u$ is an $L^p$-solution of \eqref{eq:obstacle} with $I+1$ obstacles. But, for all $n$, $\max_{i\in \{1,\dots, I\}} g^i$ is a viscosity subsolution of \eqref{eq:Dirichletinduction} and $u^n$ is a $L^p$-solution of $\eqref{eq:Dirichletinduction}$. So, the comparison principle (Proposition \ref{prop:comparison}) implies that $u^n \geq \max_{i\in \{1,\dots, I\}} g^i$ on $\bar{\mathcal{Y}}$ for all $n \in \mathbb{N}$, and, hence $u^{(I+1)} \geq \max_{i\in \{1,\dots, I\}} g^i$ on $\bar{\mathcal{Y}}$. Moreover, for all $n$, $u^n = b$ on $\partial \mathcal{Y}$, hence $u^{(I+1)} = b$ on $\partial \mathcal{Y}$. Hence, we are done if we can show that $u^{(I+1)}$ is no smaller than $g^{I+1}$ and solves $u^{(I+1)}_t + F(t,x,u^{(I+1)}, u^{(I+1)}_x, u^{(I+1)}_{xx}) = 0$ a.e. on $\{u^{(I+1)} >\max_{i\in \{1,\dots, I, I+1\}} g^i$. Both statement follows from the same steps as in the case $I=1$, hence are omitted. Therefore, $u^{(I+1)}$ is an $L^p$-solution of \eqref{eq:obstacle} with $I+1$ obstacles.

To conclude, there only remains to show that the bound \eqref{eq:W12pestimatesmoothdomainfiniteA} on the $W^{1,2,p}$-norm of $u^{(I+1)}$ holds. To this end, we first show that $u^{(I+1)}_t = g_t$, $u^{(I+1)}_x = g_x$, and $u^{(I+1)}_{xx} = g_{xx}$ a.e. on $\left\{u^{(I+1)}=g\right\}$. To see this,\footnote{The argument mimics the proof of Corollary 3.1.2.1 in \cite{evans2018measure}.} observe first that $\left\{(t,x) \in \mathcal{Y} \, :\, u^{(I+1)}(t,x) = g(t,x)\right\} = \bigcup_{a \in \{1, \dots, I\}} \mathcal{Z}^{a}$, where $\mathcal{Z}^{a} = \left\{(t,x) \in \mathcal{Y} \, :\, u^{(I+1)}(t,x) = g^a(t,x)\right\}$ and consider $\mathcal{Z}^a$, for some $a \in \{1,\dots, I\}$. By Lebesgue's theorem applied to $\mathbbm{1}_{\{ (t,x) \in \mathcal{Z}^a \}}$, for almost every $(t,x) \in \mathcal{Z}^a$,
\begin{equation}\label{eq:volume1}
    \underset{\rho\to 0^+}{\lim} \frac{B_\rho(t,x) \cap \mathcal{Z}}{B_\rho(t,x)} = 1.
\end{equation}
Moreover, $u^{(I+1)} - g^a \in {W}^{1,2,p}_{loc}(\mathcal{Y})$. Proposition A.1 in \cite{crandall1998remarks} then implies that, for almost every $(t,x) \in \mathcal{Z}^a$ and, hence, for a.e. $(t,x)$ such that \eqref{eq:volume1} holds, as $(t',x') \to (t,x)$,
\begin{equation}\label{eq:expansion}
    \begin{aligned}
    (u^{(I+1)}-g^a)(t',x') & = (u^{(I+1)}-g^a)_t (t,x)({t}'-t) + (u^{(I+1)} -g^a)_{x} \cdot ({x}' -x) \\
    & \qquad + o\left(\left| ({t}',{x}')-(t,x) \right|\right).
    \end{aligned}
\end{equation}
We next argue by contradiction that it implies that $u^{(I+1)}_t = g^a_t$ and $u^{I+1}_x = g^a_x$ on $\mathcal{Z}^a$. So, assume for a contradiction that the vector $\left((u^{(I+1)}-g^a)_t(t,x), (u^{(I+1)}-g^a)_x'\right)' = \alpha \neq 0$ at a Lebesgue point $(t,x) \in \mathcal{Z}^a$. Consider the set
\begin{equation*}
    S= \left\{ (v^t,v^x) \in \partial B_1(0,0) \, :\, \alpha \cdot v \geq \frac{1}{2} \left|\alpha\right|\right\}.
\end{equation*}
$S$ is chosen to show that $(u^{(I+1)}-g^a)>0$ on a positive measure set around $(t,x)$ and reach a contradiction. Specifically, for each $v \in S$ and $k>0$, set $(\tilde{t},\tilde{x}) = (t,x)+k v$ in \eqref{eq:expansion}. Then
\begin{equation*}
    \begin{aligned}
        (u^{(I+1)}-g)(t+k v^t, x+ kv^x) & = \alpha kv + o(\left| kv \right|) \\
        & \geq \frac{1}{2}k\left|\alpha\right| + o\left(k\right).
    \end{aligned}
\end{equation*} 
Therefore, there exists $k_0>0$ such that
\begin{equation*}
    (u^{(I+1)}-g)(t+k v^t, x+ kv^x) >0, \qquad \forall\, \, 0<k<k_0 \text{ and } v \in S, 
\end{equation*}
a contradiction to \eqref{eq:volume1}. This shows that $\alpha =0$. Hence, $u^{(I+1)}_t(t,x) = g^a_t(t,x)$ and $u^{(I+1)}_x(t,x) = g^a_x(t,x)$ a.e. on  $\mathcal{Z}^a$. By a similar  argument, $u^{(I+1)}_{xx} = g^a_{xx}$ a.e. on  $\left\{(t,x) \in \mathcal{Z}^a \, :\, u_x(t,x) = g^a_x(t,x)\right\}$, hence, a.e. on $\mathcal{Z}^a$. Since $a$ was arbitrary, the above holds for all $a$, and, hence, $u^{(I+1)}_t(t,x) = g_t(t,x)$ and $u^{(I+1)}_x(t,x) = g_x(t,x)$, and $u^{(I+1)}_{xx} = g_{xx}$ a.e. on  $\bigcup_a \mathcal{Z}^a = \left\{(t,x) \in \mathcal{Y} \, :\, u^{(I+1)}(t,x) = g(t,x)\right\}$.

Therefore, $u^{(I+1)} \in {W}^{1,2,p}\left(\mathcal{Y}\right)$ solves 
\begin{equation}\label{eq:obstacleasDirichlet}
    u^{(I+1)}_t + F\left(t,x,u^{(I+1)}, u^{(I+1)}_x, u^{(I+1)}_{xx}\right) = \mathbbm{1}_{\{ u^{(I+1)} > \underset{i=1, \dots, I+1}{\max} \, g^i(t,x) \}} h^{\iota(t,x)}(t,x)  - h^{\iota(t,x)}(t,x),
\end{equation}
where $h^i(t,x) = - g^i_t - F\left(t,x,g^i, g^i_x, g^i_{xx} \right)$ and $\iota(t,x)$ is a measurable selection from $\arg \underset{i \in \{ 1,\dots I+1 \}}{\max} \, g^{i}(t,x)$. The result then follows from Theorem 15.1.3 in \cite{krylov2018sobolev} (again using Lemma 12.1.9 in \cite{krylov2018sobolev} to bound $\left\| u \right\|_{L^{\infty}\left(\mathcal{Y}\right)}$). Thus, we have shown that when $g = \underset{a \in \left\{ 1,\dots, I, I+1 \right\}}{\max }\, g^a$ with $g^a \in W^{1,2,p}\left(\mathcal{Y}\right)$ for all $a \in \left\{ 1,\dots, I, I+1 \right\}$ and the operator $F$, domain $\mathcal{Y}$, and boundary data $b$ in \eqref{eq:obstacle} satisfy Assumptions \ref{assumption:smoothdomain}, \ref{assumption:bregularity}, and  \ref{assumption:gregularity} --\ref{assumption:growthconditions} for some $p \in (d+2, \infty)$, \eqref{eq:obstacle} has a $L^p$-solution $u \in W^{1,2,p}\left(\mathcal{Y}\right)$. Moreover, $u$ satisfies the estimate \eqref{eq:W12pestimatesmoothdomainfiniteA}.

By induction, the result holds for all $I \in \mathbb{N}$. This concludes the proof of Lemma \ref{theorem:existenceobstacleproblemfiniteA}. 	\hfill $\square$
\endproof

\vspace{1em}
Applying Theorem 12.1.7 and Lemma 12.1.9 in \cite{krylov2018sobolev} to \eqref{eq:obstacleasDirichlet}, we also obtain the following interior $W^{1,2,p}$-estimates and $L^{\infty}$ bound.
\begin{corollary}\label{corollary:interiorestimatesfiniteA}
    Suppose that Assumptions \ref{assumption:smoothdomain}, and \ref{assumption:bcontinuity}--\ref{assumption:growthconditions} hold  for some $p\in(d+2, \infty)$ and that $\left|A\right| < \infty$. If $u$ is an $L^p$-solution of \eqref{eq:obstacle}, then, for all compact subsets $\mathcal{Y}'$ of $\mathcal{Y}$,
    \begin{equation*}
        \left\| u \right\|_{W^{1,2,p}\left(\mathcal{Y'}\right)} \leq C \left( 1 + \underset{a\in A}{\sup}\, \left\| g^a \right\|_{W^{1,2,p} \left( \mathcal{Y}\right)} + \left\| G \right\|_{L^{p}\left( \mathcal{Y}\right)} + \frac{1}{dist(\mathcal{Y}', \mathcal{Y})} \left\|u \right\|_{L^{\infty}\left( \mathcal{Y} \right)} \right).
    \end{equation*}
    where $C = C\left( d, p, \lambda, \Lambda, R, diam\left(\mathcal{X}\right), T \right) \in \mathbb{R}_+$.

    Moreover, there exists $C^{\infty} = C^{\infty}\left(d, p, \lambda, \Lambda, R, \eta_F, T, diam\left(\mathcal{X}\right)\right) \in \mathbb{R}_+$ such that
	\begin{equation*}
        \left\| u \right\|_{L^{\infty} \left( \mathcal{Y}\right)} \leq C^{\infty} \left( 1 + \left\|G\right\|_{L^{p}\left( \mathcal{Y}\right)} + \underset{a\in A}{\sup}\, \left\| g^a \right\|_{W^{1,2,p} \left( \mathcal{Y}\right)} + \left\| b \right\|_{L^{\infty}\left( \mathcal{Y}\right)}  \right).
	\end{equation*}
\end{corollary}

\subsection{Proof of Theorem \ref{theorem:existenceobstacleproblemirregulardomain}}\label{sec:proofirregulardomains}

We first establish our main result for the special case of smooth domains and regular boundary functions $b$ (Section 6.1) and then prove it for the general case (Section 6.2). 

\subsubsection{Smooth Domains and Regular Boundary Functions}\label{sec:proofmaintheorem}

\begin{lemma}\label{theorem:existenceobstacleproblem}
Suppose that Assumptions \ref{assumption:smoothdomain}, \ref{assumption:bregularity}, and \ref{assumption:gregularity}--\ref{assumption:growthconditions} hold for some $p \in (d+2, \infty)$. Then, \eqref{eq:obstacle} has an $L^p$-solution $u$. Moreover, $u \in W^{1,2,p}\left(\mathcal{Y}\right)$ and there exists $C = C \left( d, p, \lambda, \Lambda, R, \eta_F, diam(\mathcal{X}), T, L^{1,Lip}(\mathcal{X}) \right) \in \mathbb{R}_+$ such that
	\begin{equation}\label{eq:W12pestimatesmoothdomain}
		\left\| u \right\|_{W^{1,2,p}\left(\mathcal{Y}\right)} \leq C \left( 1 + \underset{a \in A}{\sup} \left\| g^a \right\|_{W^{1,2,p} \left( \mathcal{Y}\right)} + \left\| b \right\|_{W^{1,2,p}\left(\mathcal{Y}\right)} + \left\| G \right\|_{L^p\left(\mathcal{Y}\right)} \right).
	\end{equation}

If, in addition, Assumption \ref{assumption:strictmonotonicity} holds, then $u$ is the unique $L^p$-solution of \eqref{eq:obstacle}.
\end{lemma}

The proof of this result builds  on Lemma~\ref{theorem:existenceobstacleproblemfiniteA}. We approximate the obstacle $g=\underset{a \in A}{\sup g^a}$ by a sequence of obstacles that satisfy the Assumptions of Lemma~\ref{theorem:existenceobstacleproblemfiniteA}. We then invoke the stability result for obstacle problems derived in Section \ref{subsec:stability} (Theorem \ref{theorem:continuityequation}) to conclude.

\vspace{0.3em}
\proof{Proof.}
    Uniqueness follows from Corollary \ref{corollary:uniqueness}. So, we only need to show existence.

    Since $A$ is separable and and $g^{\cdot}$ is continuous in $a$, there exists a countable dense subset ${A}^0 \subseteq A$ such that $\underset{a \in A}{\sup }\, g^a = \underset{a \in {A^0}}{\sup }\, g^a$ on $\mathcal{Y}$. Moreover, there exists a sequence $\left(A^{0,n}\right)_{n \in \mathbb{N}}$ of finite subsets of $A^0$ such that $\underset{a \in {A^{0,n}}}{\sup }\, g^a$ converges pointwise from below to $\underset{a \in {A^0}}{\sup }\, g^a$. 
    
    For all $n \in \mathbb{N}$, Lemma~\ref{theorem:existenceobstacleproblemfiniteA} guarantees that there exists a solution $u^n \in W^{1,2,p}\left(\mathcal{Y}\right)$ of 
    \begin{equation}\label{eq:finiteobstacles}
        \begin{cases}
            \max \left\{ u_t + F\left(t,x,u, u_x, u_{xx}\right), \,\underset{a \in A^{0,n}}{\sup}\, g^a - u \right\} = 0 \text{ on } \mathcal{Y}, \\
            u=b \text{ on } \partial \mathcal{Y}.
        \end{cases}
    \end{equation}
    Moreover, for all $n \in \mathbb{N}$,
    \begin{equation*}
		\left\| u^n \right\|_{W^{1,2,p}\left(\mathcal{Y}\right)} \leq C \left( 1 + \underset{a\in A}{\sup}\left\| g^a \right\|_{W^{1,2,p} \left( \mathcal{Y}\right)} + \left\| b \right\|_{W^{1,2,p}\left(\mathcal{Y}\right)} + \left\| G \right\|_{L^p\left(\mathcal{Y}\right)} \right).
	\end{equation*}
    Since $W^{1,2,p}(\mathcal{Y})$ is separable and reflexive, its closed bounded subsets are weakly sequentially compact by Theorem 1.32 in \cite{demengel2012functional}. Moreover $W^{1,2,p}\left(\mathcal{Y}\right)$ is compactly embedded in $\mathcal{C}^{0}\left(\bar{\mathcal{Y}}\right)$ by the Rellich-Kondrachov theorem (Theorem 2.84 in \cite{demengel2012functional}). Therefore, there exists a function $u \in W^{1,2,p}\left(\mathcal{Y}\right) \cap \mathcal{C}^{0}\left(\bar{\mathcal{Y}}\right)$ and a subsequence $\left( u^{n_j}\right)_{j \in \mathbb{N}} \subseteq \left( u^{n}\right)_{n\in \mathbb{N}}$ such that
    \begin{equation*}
        \begin{cases}
            u^{n_j}  \rightharpoonup u \text{ in } W^{1,2,p}\left( \mathcal{Y} \right),\\
            u^{n_j} \to u \text{ in } \mathcal{C}^{0}\left( \bar{\mathcal{Y}} \right),
        \end{cases}
    \end{equation*}
    as $j \to \infty$. Furthermore, $u$ satisfies the estimate \eqref{eq:W12pestimatesmoothdomain}:
    \begin{equation*}
        \begin{aligned}
            \left\| u \right\|_{W^{1,2,p}\left(\mathcal{Y}\right)} & \leq \underset{j \to \infty}{\liminf}\, \left\| u^{n_j} \right\|_{W^{1,2,p}\left(\mathcal{Y}\right)} \\
            & \leq C \left( 1 + \underset{a\in A}{\sup }\left\| g^a \right\|_{W^{1,2,p} \left( \mathcal{Y}\right)} + \left\| b \right\|_{W^{1,2,p}\left(\mathcal{Y}\right)} + \left\| G \right\|_{L^p\left(\mathcal{Y}\right)} \right).
        \end{aligned}
    \end{equation*}
    
    To conclude, there only remains to show that $u$ solves \eqref{eq:obstacle}. This follows from Theorem \ref{theorem:continuityequation}. To see this, observe that each $u^{{n_j}}$ solves the obstacle problem \eqref{eq:finiteobstacles}. Moreover, as noted, the subsequence $\left( u^{n_j}\right)_{j \in \mathbb{N}}$ converges uniformly on compact subsets of $\bar{\mathcal{Y}}$ and weakly in $W^{1,2,p}_{loc}(\mathcal{Y})$ to $u \in W^{1,2,p}_{loc}(\mathcal{Y}) \cap \mathcal{C}^0(\bar{\mathcal{Y}})$. As a result, passing to a further subsequence if necessary, we can assume that $u^{n_j}_{x}$ also converges to $u_x$ uniformly on compact subsets ${C}_D(t_0,x_0)$ by a Morrey-Sobolev embedding type theorem (Lemma 3.3, page 80 in \cite{ladyzenskaja1968linear}). Hence, by the structure condition \eqref{eq:SC}, for all fixed $v_{xx}$, ${F}\left(t, x, u^{n_j}, u^{n_j}_x, v_{xx} \right) \to F\left(t,x, u, u_x, v_{xx} \right)$ pointwise a.e.. Finally, $\underset{a \in A^{0,n}}{\sup}\, g^a \uparrow g$ pointwise. Therefore, Theorem \ref{theorem:continuityequation} applies, and $u$ solves \eqref{eq:obstacle}. \hfill $\square$
\endproof

\vspace{1em}
From Corollary \ref{corollary:interiorestimatesfiniteA} and the proof of Lemma~\ref{theorem:existenceobstacleproblem}, we obtain the following interior $W^{1,2,p}$-estimates and $L^{\infty}$ bound.
\begin{corollary}\label{corollary:interiorestimates}
    Suppose that Assumptions \ref{assumption:smoothdomain}, and \ref{assumption:bcontinuity}--\ref{assumption:strictmonotonicity} hold for some $p \in (d+2, \infty)$. If $u$ is an $L^p$-solution of \eqref{eq:obstacle}, then, for all compact subset $\mathcal{Y}'$ of $\mathcal{Y}$,
    \begin{equation*}
        \left\| u \right\|_{W^{1,2,p}\left(\mathcal{Y'}\right)} \leq C \left( 1 + \underset{a\in A}{\sup}\, \left\| g^a \right\|_{W^{1,2,p} \left( \mathcal{Y}\right)} + \left\| G \right\|_{L^{p}\left( \mathcal{Y}\right)} + \frac{1}{dist(\mathcal{Y}', \mathcal{Y})} \left\|u \right\|_{L^{\infty}\left( \mathcal{Y} \right)} \right).
    \end{equation*}
    where $C = C\left( d, p, \lambda, \Lambda, R, \eta_F, diam\left(\mathcal{X}\right), T \right) \in \mathbb{R}_+$.

    Moreover, there exists $C^{\infty} = C^{\infty}\left(d, p, \lambda, \Lambda, R, T, diam\left(\mathcal{X}\right)\right)\in \mathbb{R}_+$ such that
	\begin{equation*}
        \left\| u \right\|_{L^{\infty} \left( \mathcal{Y}\right)} \leq C^{\infty} \left( 1 + \left\|G\right\|_{L^{p}\left( \mathcal{Y}\right)} + \underset{a\in A}{\sup}\, \left\| g^a \right\|_{W^{1,2,p} \left( \mathcal{Y}\right)} + \left\| b \right\|_{L^{\infty}\left( \mathcal{Y}\right)}  \right).
	\end{equation*}
\end{corollary}

\vspace{0.3em}
\proof{Proof.}
    By Corollary \ref{corollary:uniqueness}, $u$ is the unique solution of \eqref{eq:obstacle}. Following the proof of Lemma \ref{theorem:existenceobstacleproblem}, we note that $u$ is the limit of a sequence of $W^{1,2,p}$ functions, $(u^n)_{n \in \mathbb{N}}$, such that, for all $n \in \mathbb{N}$, $u^n$ solves
    \begin{equation*}
        \begin{cases}
            \max \left\{ u^n_t + F\left(t,x,u^n, u^n_x, u^n_{xx}\right), \,\underset{a \in A^{n}}{\sup}\, g^a - u \right\} = 0 \text{ on } \mathcal{Y}, \\
            u^n=b^n \text{ on } \partial \mathcal{Y},
        \end{cases}
    \end{equation*}
    where $A^n$ is finite and the $W^{1,2,p}$ functions $b^n$, $n\in \mathbb{N}$, converge to $b$. Corollary \ref{corollary:interiorestimatesfiniteA} then implies that, for all $n \in \mathbb{N}$, 
    \begin{equation*}
        \left\| u^n \right\|_{L^{\infty} \left( \mathcal{Y}\right)} \leq C^{\infty} \left( 1 + \left\|G\right\|_{L^{p}\left( \mathcal{Y}\right)} + \underset{a\in A}{\sup}\, \left\| g^a \right\|_{W^{1,2,p} \left( \mathcal{Y}\right)} + \left\| b \right\|_{L^{\infty}\left( \mathcal{Y}\right)}  \right)
	\end{equation*}
    for some $C^{\infty} = C^{\infty}\left(d, p, \lambda, \Lambda, R, T, diam\left(\mathcal{X}\right)\right) \in \mathbb{R}_+$. Moreover, for all $n \in \mathbb{N}$ and all compact subset $\mathcal{Y}'$ of $\mathcal{Y}$,
    \begin{equation*}
        \left\| u^n \right\|_{W^{1,2,p}\left(\mathcal{Y'}\right)} \leq C \left( 1 + \underset{a\in A}{\sup}\, \left\| g^a \right\|_{W^{1,2,p} \left( \mathcal{Y}\right)} + \left\| G \right\|_{L^{p}\left( \mathcal{Y}\right)} + \frac{1}{dist(\mathcal{Y}', \mathcal{Y})} \left\|u^n \right\|_{L^{\infty}\left( \mathcal{Y} \right)} \right),
    \end{equation*}
    for some $C = C\left( d, p, \lambda, \Lambda, R, \eta_F, diam\left(\mathcal{X}\right), T \right) \in \mathbb{R}_+$.
	
    The desired results then follow, since
    \begin{equation*}
          \left\| u \right\|_{L^{\infty}\left(\mathcal{Y}\right)} \leq \underset{n \to \infty}{\liminf}\, \left\| u^{n} \right\|_{L^{\infty}\left(\mathcal{Y}\right)} \text{ and } \left\| u \right\|_{W^{1,2,p}\left(\mathcal{Y'}\right)} \leq \underset{n \to \infty}{\liminf}\, \left\| u^{n} \right\|_{W^{1,2,p}\left(\mathcal{Y'}\right)}.
    \end{equation*} \hfill $\square$
\endproof

\subsubsection{Proof of Theorem~\ref{theorem:existenceobstacleproblemirregulardomain}}

Lemma~\ref{theorem:existenceobstacleproblem} guarantees that the obstacle problem has an $L^p$-solution when $\mathcal{X}$ and $b$ satisfy Assumption \ref{assumption:smoothdomain} and \ref{assumption:bregularity}. To generalize the result to the weaker Assumptions~\ref{assumption:lipschitzdomain} and \ref{assumption:bcontinuity}, we study a sequence of equations, each satisfying the assumptions of Lemma~\ref{theorem:existenceobstacleproblem}, that converges to the equation \eqref{eq:obstacle}. In particular, we approximate (i) $\mathcal{X}$ by a sequence of smooth domains whose cone parameters are uniformly controlled, and (ii) $b$ by a sequence of equicontinuous functions in $W^{1,2,p}\left( \mathcal{Y} \right)$. The $L^p$ solutions of the equations in the approximating sequence form an equicontinuous (by Lemma \ref{lemma:localmodulusofcontinuity}, below) and weakly compact (in $W^{1,2,p}_{loc}\left(\mathcal{Y}\right)$) family. The Arzel\`{a}-Ascoli theorem then guarantees that a subsequence converges uniformly on the compact subset of $\bar{\mathcal{Y}}$ to some function in $W^{1,2,p}_{loc}\left(\mathcal{Y}\right) \cap \mathcal{C}^0\left(\Bar{\mathcal{Y}}\right)$. Finally, we invoke the stability result for obstacle problems (Theorem \ref{theorem:continuityequation}) derived in Section \ref{subsec:stability} to conclude.

\begin{lemma}\label{lemma:localmodulusofcontinuity}
    Suppose that Assumptions \ref{assumption:smoothdomain}, \ref{assumption:bregularity}, and \ref{assumption:gregularity}--\ref{assumption:strictmonotonicity} hold for some $p \in (d+2, \infty)$. Let $u$ be an $L^p$-solution of \eqref{eq:obstacle}. For all $D>0$, there exists a family of modulus of continuity $\bar{\omega}_D$ such that, for all $(t,x), (t',x') \in \bar{\mathcal{Y}} \cap C_D(0,0)$,
    \begin{equation}\label{eq:modulusofcontinuity}
        \left|u(t,x) - u(t',x')\right| \leq \bar{\omega}_{D}\left(\left| (t,x) - (t',x') \right|\right),
    \end{equation}
    where $\bar{\omega}_D$ depends only on $p,d, \lambda, \Lambda, R, \eta_F, D$, $\underset{a \in A}{\sup \, }\left\| g^a\right\|_{W^{1,2,p}\left(\mathcal{Y}\cap C_{D+1}(0,0)\right)}$, $\left\| G\right\|_{L^{p}\left(\mathcal{Y}\cap C_{D+1}(0,0)\right)}$, $\left\|u\right\|_{L^{\infty}\left( \mathcal{Y} \cap C_{D+1}(0,0) \right)}$, the modulus of continuity of $b$ on $\partial \mathcal{Y} \cap C_{D+1}(0,0)$, and the parameters of the cone condition of $\mathcal{X}$ on $ \mathcal{X} \cap B_{D+1}(0)$. 
\end{lemma}

\begin{remark}\label{remark:localmodulusofcontinuity}
    The modulus of continuity $\bar{\omega}_D$ for the solution of \eqref{eq:obstacle} on $\bar{\mathcal{Y}} \cap C_D(0,0)$ given in Lemma \ref{lemma:localmodulusofcontinuity} only depends on the characteristics of the primitive inside the cylinder $C_{D+1}(0,0)$. It is independent of $\underset{a \in A}{\sup \, }\left\| g^a\right\|_{W^{1,2,p}\left(\mathcal{Y} \setminus C_{D+1}(0,0)\right)}$, $\left\| G\right\|_{L^{p}\left(\mathcal{Y} \setminus C_{D+1}(0,0)\right)}$, $\left\|u\right\|_{L^{\infty}\left( \mathcal{Y} \setminus C_{D+1}(0,0) \right)}$, the modulus of continuity of $b$ on $\partial \mathcal{Y} \setminus C_{D+1}(0,0)$, and of the parameters of the cone condition of $\mathcal{X}$ on $ \mathcal{X} \setminus B_{D+1}(0)$.
\end{remark}

\vspace{0.3em}
\proof{Proof.}
By the structure condition \eqref{eq:SC}, for any $L^p$-solution $u$ of \eqref{eq:obstacle} and almost all $(t,x) \in \mathcal{Y}\cap C_{D+1}(0,0)$, we have
\begin{equation*}
    \begin{aligned}
        \max& \left\{u_t + F(t,x, u, u_x, u_{xx}), \, g-u\right\} \\
        & \leq u_t + \mathcal{P}^+_{\lambda, \Lambda}\left( u_{xx} \right)  + R  \left| u_x\right| \\
        & \qquad + \bar{C} \left(1+ \left\|u\right\|_{L^{\infty}\left(\mathcal{Y} \cap C_{D+1}(0,0)\right)}+ \sup_{a \in A} \left( \left| g^a\right| + \left| g^a_t\right| + \left| g^a_x\right| + \left| g^a_{xx}\right| \right) + G(t,x) \right),
    \end{aligned}
\end{equation*}
for some $\bar{C} > 0$ depending only on $p,d,\lambda, \Lambda, R$, and $D$. Similarly,
\begin{equation*}
\begin{aligned}
    \max&\left\{u_t + F(t,x, u, u_x, u_{xx}), \, g-u\right\} \\
        & \geq u_t + \mathcal{P}^-_{\lambda, \Lambda}\left( u_{xx} \right)  - R  \left| u_x\right|  - \bar{C} \left(1+ \left\|u\right\|_{L^{\infty}\left(\mathcal{Y}\right)} + \sup_{a \in A} \left( \left| g^a\right| + \left| g^a_t\right| + \left| g^a_x\right| + \left| g^a_{xx}\right| \right) + G(t,x) \right).
\end{aligned}
\end{equation*}
    Next, we define two boundary values $\bar{u}$ and $\underline{u}$ on $\partial \left( \mathcal{Y} \cap C_{D+1}(0,0)\right)$ by
    \begin{equation*}
        \bar{u}(t,x) = \begin{cases}
            b(t,x) \text{ if } (t,x) \in \left(\partial \mathcal{Y} \right) \cap C_D(0,0), \\
            dist\left(\{(t,x)\}, C_{D+1}\right) b(t,x) + \left(1- dist\left(\{(t,x)\}, C_{D+1}\right) \right) \left\| u \right\|_{L^{\infty}\left(\mathcal{Y} \cap C_{D+1}(0,0)\right)} \\
            \qquad \qquad \text{ if } (t,x) \in \partial \mathcal{Y} \cap \left(C_{D+1} \setminus C_D(0,0)\right) \\
            \left\| u \right\|_{L^{\infty}\left(\mathcal{Y} \cap C_{D+1}\right)} \text{ if } (t,x) \in \partial \left( \mathcal{Y} \cap  C_{D+1}\right) \setminus \partial \mathcal{Y}.
        \end{cases}
    \end{equation*}
    and
    \begin{equation*}
        \underline{u}(t,x) = \begin{cases}
            b(t,x) \text{ if } (t,x) \in \left(\partial \mathcal{Y} \right) \cap C_D(0,0), \\
            dist\left(\{(t,x)\}, C_{D+1}\right) b(t,x) - \left(1- dist\left(\{(t,x)\}, C_{D+1}\right) \right) \left\| u \right\|_{L^{\infty}\left(\mathcal{Y} \cap C_{D+1}(0,0)\right)} \\
            \qquad \qquad \text{ if } (t,x) \in \partial \mathcal{Y} \cap \left(C_{D+1} \setminus C_D(0,0)\right) \\
            - \left\| u \right\|_{L^{\infty}\left(\mathcal{Y} \cap C_{D+1}(0,0)\right)}  (t,x) \in \partial \left( \mathcal{Y} \cap  C_{D+1}(0,0)\right) \setminus \partial \mathcal{Y}.
        \end{cases}
    \end{equation*}
    Both $\bar{u}$ and $\underline{u}$ are equal to $b$ on $\partial \mathcal{Y} \cap C_D(0,0)$ and equal to a convex combination of $b$ and the constants $\left\| u \right\|_{L^{\infty}(\mathcal{Y} \cap C_{D+1})}$ and $-\left\| u \right\|_{L^{\infty}(\mathcal{Y} \cap C_{D+1})}$, respectively, on $\partial \mathcal{Y} \cap \left( \bar{C}_{D+1} \setminus C_D\right)$. In particular, we defined $\bar{u}$ and $\underline{u}$ so that (i) they are continuous and their moduli of continuity are controlled by $\left\| u \right\|_{L^{\infty}\left(\mathcal{Y} \cap C_{D+1}(0,0)\right)}$ and the modulus of continuity of $b$ on $\partial Y \cap C_{D+1}(0,0)$, and (ii) $\bar{u} \geq u \geq \underline{u}$ on $\partial \left( \mathcal{Y} \cap C_{D+1}(0,0)\right)$.
    
    It follows that $u$ is also an $L^p$-subsolution of
	\begin{equation}\label{eq:boundarycontinuitysubsolution}
		\begin{cases}
			& v_t + \mathcal{P}^+_{\lambda, \Lambda}\left( v_{xx} \right)  + R  \left| v_x\right| \\  
            & \qquad = - \bar{C} \left(1+ \left\|u\right\|_{L^{\infty}\left( \mathcal{Y} \cap C_{D+1}(0,0)\right)} + \sup_{a \in A} \left( \left| g^a\right| + \left| g^a_t\right| + \left| g^a_x\right| + \left| g^a_{xx}\right| \right) + G(t,x) \right) \\
            & \qquad \qquad \text{ on } \mathcal{Y} \cap C_{D+1}(0,0), \\
			& v = \bar{u} \text{ on } \partial \mathcal{Y} \cap C_{D+1}(0,0),
		\end{cases}
	\end{equation}
	and an $L^p$-supersolution of
	\begin{equation}\label{eq:boundarycontinuitysupersolution}
		\begin{cases}
			& v_t + \mathcal{P}^-_{\lambda, \Lambda}\left( v_{xx} \right)  - R  \left| v_x\right| \\  
            & \qquad = \bar{C} \left(1+ \left\|u\right\|_{L^{\infty}\left(\mathcal{Y}\cap C_{D+1}(0,0)\right)} + \sup_{a \in A} \left( \left| g^a\right| + \left| g^a_t\right| + \left| g^a_x\right| + \left| g^a_{xx}\right| \right) + G(t,x) \right) \\
            & \qquad \qquad \text{ on } \mathcal{Y}\cap C_{D+1}(0,0), \\
			& v = \underline{u} \text{ on } \partial \mathcal{Y}\cap C_{D+1}(0,0).
		\end{cases}
	\end{equation}
    By Theorem 4.5 in \cite{crandall1999existence}, there exists a viscosity solution $\bar{U}$ of \eqref{eq:boundarycontinuitysubsolution} and a viscosity solution $\underline{U}$ of \eqref{eq:boundarycontinuitysupersolution}. Moreover, by Proposition 4.6 and Remark 4.3 in \cite{crandall1999existence}, there exists a modulus of continuity $\bar{\omega}_{boundary}$ depending only on $p,d,\lambda, \Lambda, R, D, \underset{a \in A}{\sup \, }\left\| g^a\right\|_{W^{1,2,p}\left(\mathcal{Y} \cap C_{D+1}(0,0)\right)}$, $\left\| G\right\|_{L^{p}\left(\mathcal{Y}\cap C_{D+1}(0,0)\right)}, \left\|u\right\|_{L^{\infty}\left(\mathcal{Y} \cap C_{D+1}(0,0) \right)}$, the modulus of continuity of $b$ on $\mathcal{Y} \cap C_{D+1}(0,0)$, and the parameters of the cone condition of $\mathcal{X}$ on $\mathcal{X} \cap B_{D+1}(0)$ such that
    \begin{align*}
        \max\left\{ \left| \bar{U}(t,x) - b\left(\bar{t}, \bar{x}\right)\right|,  \left| \underline{U}(t,x) - b\left(\bar{t}, \bar{x}\right)\right| \right\} \leq \bar{\omega}_{boundary}\left(\left|(t,x) - (\bar{t}, \bar{x}) \right| \right).
    \end{align*}
    for all $(t,x) \in \bar{\mathcal{Y}}\cap C_{D}(0,0)$ and all $\left(\bar{t}, \bar{x}\right) \partial \mathcal{Y} \cap C_{D}(0,0)$. By our comparison principle (Proposition \ref{prop:comparison}), since $\underline{u} \leq u \leq \bar{u}$ on $\partial \mathcal{Y}\cap C_D(0,0)$,
	\begin{equation*}
		\underline{U} \leq u \leq \bar{U} \text{ on } \bar{\mathcal{Y}} \cap C_{D+1}(0,0).
	\end{equation*}
	Therefore, for all $(t,x) \in \bar{\mathcal{Y}}\cap C_{D}(0,0)$ and all $\left(\bar{t}, \bar{x}\right) \partial \mathcal{Y} \cap C_{D}(0,0)$,
    \begin{equation}\label{eq:modulusofcontinuityattheboundary}
	   \left| u(t,x) - b\left(\bar{t}, \bar{x}\right)\right| \leq \max\left\{ \left| \bar{U}(t,x) - b\left(\bar{t}, \bar{x}\right)\right|,  \left| \underline{U}(t,x) - b\left(\bar{t}, \bar{x}\right)\right| \right\} \leq \bar{\omega}_{boundary}\left(\left|(t,x) - (\bar{t}, \bar{x}) \right| \right)
	\end{equation}
    where $\bar{\omega}_{boundary}$ is the modulus of continuity defined above.
    
    Finally, by Corollary \ref{corollary:interiorestimates} and a Morrey-Sobolev embedding theorem (Theorem 2.84 in \cite{demengel2012functional}), $u$ is $\alpha$-H\"{o}lder continuous on any compact subset $\mathcal{Y}'$ of $\mathcal{Y} \cap C_{D+1}\left(0,0\right)$ for all $\alpha < 1 - (d+1/p)$, with H\"{o}lder norm bounded above by 
    \begin{equation*}
        C^{H} \left( 1 +  \underset{a\in A}{\sup }\left\| g^a \right\|_{W^{1,2,p} \left( \mathcal{Y}\cap C_{D+1}\left(0,0\right)\right)} + \left\| u \right\|_{L^{\infty}\left(\mathcal{Y} \cap C_{D+1}\left(0,0\right)\right)} + \left\| G \right\|_{L^p\left(\mathcal{Y}\cap C_{D+1}\left(0,0\right)\right)}  \right)
    \end{equation*}
    where $C^H$ depends only, as above, on $p,d,\lambda, \Lambda, R, \eta_F, D$, and on the distance of the compact set to the boundary $dist(\mathcal{Y}', \partial\mathcal{Y})$. 
    
    Combining these two estimates, we obtain the desired modulus of continuity $\bar{\omega}_D$. \hfill $\square$ \endproof

\vspace{1em}
{\bf Proof of Theorem \ref{theorem:existenceobstacleproblemirregulardomain}.}

Uniqueness follows from Corollary \ref{corollary:uniqueness}. The remainder of the proof shows existence and establishes the interior estimates \eqref{eq:W12pestimatesirregulardomain}.
    
By Theorem 5.1 in \cite{doktor1976approximation}, there exists a strictly increasing (in the sense of set inclusion) sequence of bounded smooth open subsets of $\mathcal{X}$, $\left( \mathcal{X}_n\right)_{n \in \mathbb{N}}$, and a constant $M>0$ such that, for all $n \in \mathbb{N}$,
\begin{equation*}
        \max \left\{ L^{0,Lip}(\mathcal{X}), \, L^{0,Lip}(\mathcal{X}_n) \right\} \leq M.
\end{equation*}
By Theorem 1.2.2.2 in \cite{grisvard2011elliptic}, $\mathcal{X}$ and $\mathcal{X}_n$, $n \in \mathbb{N}$ satisfy a uniform exterior cone condition of size $(h, \theta)$, for some $h,\theta >0$ independent of $n$. In particular $(h, \theta)$ depends only on $M$ and $diam(\mathcal{X})$.
    
Define $\mathcal{Y}_n =  \left[0, \frac{2n-1}{2n}T \right) \times \mathcal{X}_n$, and consider the following sequence of nonlinear obstacle problems: for all $n \in \mathbb{N}$, 
\begin{equation}\label{eq:Eqnirregulardomain}
		\begin{cases}
		  \max \left\{ u_t+ F\left(t, x, u, u_x, u_{xx} \right), \, g-u\right\} = 0 \text{ on } \mathcal{Y}_n \\
		  u = b^n \text{ on } \partial \mathcal{Y}_n,
		\end{cases} 
\end{equation}
where $b_n \in W^{1,2,p}(\mathcal{Y})$ for all $n\in \mathbb{N}$, and $b^n \to b \in \mathcal{C}^0(\bar{Y})$.

By Lemma~\ref{theorem:existenceobstacleproblem}, for all $n \in \mathbb{N}$, there exists a unique $L^p$-solution $u^n \in W^{1,2,p}\left(\mathcal{Y}_n\right) \subset W^{1,2,p}_{loc}\left( \mathcal{Y}_n \right) \cap \mathcal{C}^0\left( \bar{\mathcal{Y}}_n \right)$ of \eqref{eq:Eqnirregulardomain}. Extending the functions $u^n$, $n\in \mathbb{N}$, by letting $u^n(t,x) = b(t,x)$ on $\bar{\mathcal{Y}} \setminus \mathcal{Y}_n$, we obtain a sequence $\left( u^n\right)_{n \in \mathbb{N}} \subset \mathcal{C}^0\left( \bar{\mathcal{Y}} \right)$.
	
Next we show that the sequence $\left( u^n \right)_{n \in \mathbb{N}}$ has a subsequence that converges (i) uniformly in $\mathcal{C}^0\left( \bar{\mathcal{Y}} \right)$ and (ii) weakly in ${W}^{1,2,p}_{loc}\left( \mathcal{Y} \right)$ to some function $u \in {W}^{1,2,p}_{loc}\left( \mathcal{Y} \right) \cap \mathcal{C}^0\left(\bar{\mathcal{Y}}\right)$. 

We start with (i). By Lemma \ref{lemma:localmodulusofcontinuity} with $D> T\vee diam(\mathcal{X})$ (using Lemma 12.1.9 in \cite{krylov2018sobolev} to uniformly control $\left\|u^n\right\|_{L^{\infty} \left(\mathcal{Y}\right)}$, $n \in \mathbb{N}$), the sequence $\left( u^{n}\right)_{n \in \mathbb{N}}$ is equicontinuous. As a result, it has a convergent subsequence in $\mathcal{C}^0\left(\bar{\mathcal{Y}}\right)$ by the Arzel\`{a}-Ascoli theorem (Theorem A.5 in \cite{rudin1973functional}), which, with a small abuse of notation, we relabel $\left( u^n\right)_{n \in \mathbb{N}}$.

Next, we show (ii), i.e., that the subsequence above, which converges in $\mathcal{C}^0\left(\bar{\mathcal{Y}}\right)$, has a weakly convergent subsequence in ${W}^{1,2,p}_{loc}\left( \mathcal{Y} \right)$. To do so, we prove that the $W^{1,2,p}$-norms of the functions $u^n$, $n\in \mathbb{N}$, restricted to any given compact subset of $\mathcal{Y}$, are eventually uniformly bounded. 
    
It is easy to see that, for all $n \geq 3$, $u^n \in W^{1,2,p}_{loc}\left( \mathcal{Y}_n\right) \cap \mathcal{C}^0\left( \bar{\mathcal{Y}}_n \right)$ is also an $L^p$-solution of
\begin{equation*}
        \begin{cases}
            \max \left\{ u_t+ F\left(t, x, u, u_x, u_{xx} \right), \, g-u\right\} = 0 \text{ on } \mathcal{Y}_{n'} \\
            u(t,x) = u^n(t,x) \text{ on } \partial \mathcal{Y}_{n'},
        \end{cases}
\end{equation*}
for all $n' \leq n-2$. So, by Corollary \ref{corollary:interiorestimates}, for all $n \geq 3$ and all $n' \leq n -2$,
\begin{equation}\label{eq:interiorestimatesn'2}
        \left\| u^n \right\|_{W^{1,2,p}\left(\mathcal{Y}_{n'}\right)} \leq C^{n'} \left( 1+ \left\|G\right\|_{L^{p}\left( \mathcal{Y}\right)} + \underset{a\in A}{\sup}\, \left\| g^a \right\|_{W^{1,2,p} \left( \mathcal{Y}\right)} + \left\| b \right\|_{L^{\infty}\left( \mathcal{Y}\right)}\right),
\end{equation}
for some $C^{n'} =  C^{n'}\left( d, p, \lambda, \Lambda, R, \eta_F, diam\left(\mathcal{X}\right), T, dist(\mathcal{Y}_{n'}, \mathcal{Y}_{n'+1}) \right) \in \mathbb{R}_+$. Crucially, the right-hand side of \eqref{eq:interiorestimatesn'2} is independent of $n$. 
    
We then show that there exists a subsequence of $\left(u^n\right)_{n \in \mathbb{N}}$ that converges to the limit $u \in \mathcal{C}^0\left( \bar{\mathcal{Y}} \right)$ obtained in (i) by a diagonal argument, and, therefore, that $u \in {W}^{1,2,p}_{loc}\left(\mathcal{Y} \right)\cap\mathcal{C}^0\left( \bar{\mathcal{Y}} \right)$. For all $n \geq 3$,\eqref{eq:interiorestimatesn'2} implies that
\begin{equation*}
        \left\| u^n \right\|_{W^{1,2,p}\left(\mathcal{Y}_{1}\right)} \leq C^{1} \left( 1+ \left\|G\right\|_{L^{p}\left( \mathcal{Y}\right)} + \underset{a\in A}{\sup}\, \left\| g^a \right\|_{W^{1,2,p} \left( \mathcal{Y}\right)} + \left\| b \right\|_{L^{\infty}\left( \mathcal{Y}\right)}\right).
\end{equation*}
${W}^{1,2,p}\left(\mathcal{Y}_1 \right)$ is separable and reflexive, and, hence, its closed bounded subsets are weakly sequentially compact by Theorem 1.32 in \cite{demengel2012functional}. Moreover $W^{1,2,p}(\mathcal{Y}_1)$ is compactly embedded in $\mathcal{C}^0 \left(\mathcal{Y}_1 \right)$ by the Rellich-Kondrachov theorem (Theorem 2.84 in \cite{demengel2012functional}). Therefore, there exists a function $\tilde{u}$ defined on $\bar{\mathcal{Y}}$ which restriction on $\mathcal{Y}_1$ is continuous, and a subsequence $\left( u^{n_j}\right)_{j \in \mathbb{N}} \subseteq \left( u^{n}\right)_{n\in \mathbb{N}}$ such that
\begin{equation*}
		\begin{cases}
			u^{n_j} \rightharpoonup \tilde{u} \text{ in } W^{1,2,p}\left( \mathcal{Y}_1 \right)\\
			u^{n_j} \to \tilde{u} \text{ in } \mathcal{C}^0\left( \bar{\mathcal{Y}}_1 \right),
		\end{cases}
\end{equation*}
as $j \to \infty$. In particular, $\tilde{u}$ must coincide with $u$ on $\mathcal{Y}_1$.
		
Proceeding with a diagonal argument, we see that there exists a subsequence $\left( u^{n_j}\right)_{j \in \mathbb{N}} \subseteq \left( u_{n}\right)_{n\in \mathbb{N}}$ such that
\begin{equation*}
		\begin{cases}
    		u^{n_j}  \rightharpoonup u \text{ in } W^{1,2,p}_{loc}\left( \mathcal{Y} \right) \\
    		u^{n_j}  \to u \text{ in } \mathcal{C}^0\left( \bar{\mathcal{Y}} \right),
		\end{cases}
\end{equation*}
as $j \to \infty$. Moreover, $u$ satisfies the interior estimates \eqref{eq:W12pestimatesirregulardomain} since, for all $n \in \mathbb{N}$,
\begin{equation*}
    \begin{aligned}
        \left\| u \right\|_{W^{1,2,p}\left(\mathcal{Y}_{n}\right)} & \leq \underset{j \to \infty}{\liminf}\, \left\| u^{n_j} \right\|_{W^{1,2,p}\left(\mathcal{Y}_{n}\right)} \\
        & \leq C^{n} \left( 1+ \left\|G\right\|_{L^{p}\left( \mathcal{Y}\right)} + \underset{a\in A}{\sup}\, \left\| g^a \right\|_{W^{1,2,p} \left( \mathcal{Y}\right)} + \left\| b \right\|_{L^{\infty}\left( \mathcal{Y}\right)}\right).
    \end{aligned}
\end{equation*}
That is, $u \in W^{1,2,p}_{loc}\left( \mathcal{Y} \right) \cap \mathcal{C}^0\left( \bar{\mathcal{Y}} \right)$ and $u$ satisfies the interior $W^{1,2,p}$-estimates \eqref{eq:W12pestimatesirregulardomain}.
		
There remains to show that $u \in W^{1,2,p}_{loc}\left( \mathcal{Y}\right) \cap \mathcal{C}^0\left( \bar{\mathcal{Y}} \right)$ is an $L^p$-solution of \eqref{eq:obstacle}. By construction, $\left( u^{n_j}\right)_{j \in \mathbb{N}} \subseteq W^{1,2,p}_{loc} \left( \mathcal{Y} \right)\cap \mathcal{C}^0\left( \bar{\mathcal{Y}} \right)$ converges to $u$ weakly in $W^{1,2,p}_{loc}\left( \mathcal{Y}\right) $ and uniformly on $\bar{\mathcal{Y}}$. The result then follows from Theorem~\ref{theorem:continuityequation}.\hfill $\square$

% Appendix here
% Options are (1) APPENDIX (with or without general title) or
%             (2) APPENDICES (if it has more than one unrelated sections)
% Outcomment the appropriate case if necessary
%
% \begin{APPENDIX}{<Title of the Appendix>}
% \end{APPENDIX}
%
%   or
%
% \begin{APPENDICES}
% \section{<Title of Section A>}
% \section{<Title of Section B>}
% etc
% \end{APPENDICES}

% References here (outcomment the appropriate case)

% CASE 1: BiBTeX used to constantly update the references
%   (while the paper is being written).
%\bibliographystyle{informs2014} % outcomment this and next line in Case 1
%\bibliography{<your bib file(s)>} % if more than one, comma separated

%\bibliographystyle{informs2014} % outcomment this and next line in Case 1
%\bibliography{sample} % if more than one, comma separated

% CASE 2: BiBTeX used to generate mypaper.bbl (to be further fine tuned)
%\input{mypaper.bbl} % outcomment this line in Case 2

%If you don't use BiBTex, you can manually itemize references as shown below.

%\bibliographystyle{nonumber}

\bibliography{biblio}

%%%%%%%%%%%%%%%%%
\end{document}